\documentclass[letterpaper, paper,11pt]{AAS}	% for preprint proceedings

\usepackage{bm}
\usepackage{amsmath}
\usepackage{svg}
\usepackage{subfigure}
\usepackage[colorlinks=true, pdfstartview=FitV, linkcolor=black, citecolor= black, urlcolor= black]{hyperref}
\usepackage{overcite}
\usepackage{footnpag}			      	% make footnote symbols restart on each page
\usepackage{arydshln}

\DeclareMathOperator*{\minimize}{minimize}

\PaperNumber{22-810}

\begin{document}

\title{Trajectory-Planning Attitude Control System for Satellites with Magnetic Attitude Control and One Reaction Wheel}

\author{Patrick McKeen\thanks{Ph.D. Candidate, Department of Aeronautics and Astronautics, 77
Massachusetts Avenue; pmckeen@mit.edu (Corresponding Author)},  
Alex Meredith\thanks{Master's Student, Department of Aeronautics and Astronautics, 77
Massachusetts Avenue; ameredit@mit.edu },
\ and Kerri Cahoy\thanks{Associate Professor, Department of Aeronautics and Astronautics, 77
Massachusetts Avenue. Associate Fellow AIAA}
}

\maketitle{}

\begin{abstract}
Selecting attitude actuators for spacecraft in low earth orbit (LEO) often involves a trade between reaction wheels and magnetorquers. Reaction wheels generally have higher size, weight and power as well as cost, compared with magnetorquers, but have performance advantages. We consider a hybrid system relying on a single reaction wheel and three magnetorquers, using a trajectory planning control method that focuses on meeting pointing goals and minimizing the need to desaturate the reaction wheel. The use of planning allows for intelligent use of the reaction wheel and accommodates time-changing goals and orbit-varying magnetic fields. Preliminary results show trajectories with fractional-degree-accuracy, as well as management of reaction wheel momentum and avoidance of stiction. We also compare to a magnetorquer-only planning-based control system and show better convergence and faster slews.
\end{abstract}
\section{Introduction}
Attitude control is a critical spacecraft capability. Early satellites often relied on methods such as spin-stabilization, gravity gradient control, or permanent magnetic control to control attitude.\cite{wertz_space_2011,cappelletti_cubesat_2020} As space technology has developed, active control has become more common, typically using thrusters, reaction wheels, or control moment gyroscopes. These control systems have all been flown functionally on many missions, but typically come with a trade-off in terms of size, weight, power, and cost (SWaP-C). In recent years, the rise of CubeSats has fueled increasing demand for innovative attitude control systems with low SWaP-C. This has led to interest in magnetic control systems, which have been used on satellites for decades, but rarely as primary control. Magnetic control systems have a specific limitation in that they are underactuated and can only be used in orbit around planetary bodies with magnetic fields. Magnetorquers cannot torque around the local magnetic field line, but their advantages are numerous: they cost less, do not require fuel, do not saturate as reaction wheels or control moment gyroscopes can, do not create vibration as spinning wheels do, and are more resistant to failure (due to lack of moving parts). They also tend to be small in size and weight, and use no power when not actuating. These advantages have spurred interest in how their underactuation can be addressed or overcome.

Magnetic control has been used on a variety of satellites, including several which primarily or exclusively used magnetic control: the TANGO satellite in the PRISMA mission \cite{chasset_3-axis_2006,chasset_3-axis_2013}, the GOCE spacecraft \cite{sechi_magnetic_2005} and the \O{}ersted satellite\cite{wisniewski_three-axis_1996,wisniewski_satellite_1997,wisniewski_linear_2000}. In addition to these magnetically-controlled missions, many potential magnetic controllers have been proposed and described in literature. A good (and regularly updated) survey can be found in Ovchinnikov \& Roldugin \cite{ovchinnikov_survey_2019}. However, many of these algorithms are limited. Many use external torques, such as gravity gradient or solar radiation pressure, to stabilize the satellite.\cite{wisniewski_optimal_1999,psiaki_magnetic_2001,stuurman_ryefemsat_2010,frik_attitude_1970,belokonov_analysis_2018}
Most of these methods require a specific satellite geometry or apply to specific orbits or pointing goals. Even those that do not\cite{lovera_global_2006} are based on orbit-average magnetic fields and thus can take a long time to stabilize.
Other methods effectively use Lyapunov or Proportional-Derivative (PD) control, adapted for the underactuation limitation. 
Lovera and Astolfi proposed a locally exponentially stable magnetic control law that projects the planned control torque onto the plane of available torque \cite{lovera_global_2005}. Other methods have been shown to be robust to uncertainty; these can be effective, but can also take a long time to complete a slew or adapt to changing conditions and goals\cite{ivanov_advanced_2017}.  Another common method is to use a sliding mode. 
Several sliding-mode magnetic controllers have been proposed, including versions demonstrated on both the TANGO and \O{}ersted spacecraft \cite{wisniewski_sliding_1998,wang_satellite_1998,chasset_3-axis_2006,janardhanan_attitude_2012,sofyali_purely_2012}.

Many of these approaches are elegant and successful, but advances in the field of robotics have created trajectory-optimization algorithms that can find and generate orientation trajectories that match goals, even in the underactuated case.\cite{liang_time-optimal_2004,kumar_h_2016,krogstad_explicit_2005}. The approach most relevant to this paper is the one proposed by Gatherer and Manchester \cite{gatherer2018magnetorquer}, which applies trajectory optimization techniques from robotics to the problem of controlling magnetorquer-only satellites. In particular, Gatherer and Manchester showed that ALTRO \cite{jackson_al_ilqr_tutorialpdf_2019,jackson_planning_2021,howell_altro_2019}, a trajectory optimization algorithm based on an augmented Lagrangian version of an iterative linear quadratic regulator (iLQR), can plan $180^\circ$ slews for a 3U satellite in an arbitrary low Earth orbit that take place in less than ten minutes. 
This use of explicit trajectory optimization is implemented on the  BeaverCube-1 (BC) Mission, a 3U satellite developed by the Space Telecommunications, Astronomy, and Radiation (STAR) Lab at MIT, which launched in July 2022, and plans to test this method on-orbit.\cite{m_campbell_beavercube_2019} Specifically, this method was used because it has the potential to deliberately plan around and counter disturbance torques (such as from propulsion).
\\ \\
While there are many interesting approaches to magnetic-only control, the development process of BC showed that there are fundamental limitations to magnetic control, due to underactuation and low torque, that render responsive and immediate 3-axis control difficult. This inspired the investigation of adding a single reaction wheel to a satellite. In this paper, we use an adapted version of ALTRO developed for BeaverCube-1 for trajectory planning, and extend this version of ALTRO to plan slews for a satellite with three magnetorquers and a single reaction wheel (1RW). We demonstrate that ALTRO can plan single-target slews for a 1RW satellite with sub-degree precision in the vast majority of cases, and can plan multi-target slews for a 1RW satellite with sub-degree precision in a significant fraction of cases, all while limiting angular momentum build-up. This represents a great improvement over the capabilities of magnetorquer-only satellites.

\section{Background}
In this paper, we use the Hamilton convention for the unit attitude quaternion \cite{sola2017quaternion}, such that if the rotation from the body-fixed satellite frame to the global ECI frame is given by a rotation of $\theta$ radians about an axis given by the unit vector $[x\ y\ z ]^T$, the unit attitude quaternion $q$ is as follows:
\begin{equation} \label{eq:ham_convention}
    \begin{bmatrix} q_w \\ q_x \\ q_y \\ q_z \end{bmatrix} = \begin{bmatrix} \cos(\frac\theta2) \\ x\sin(\frac\theta2) \\ y\sin(\frac\theta2) \\ z\sin(\frac\theta2) \end{bmatrix}
\end{equation}
The kinematics of the unit attitude quaternion $q$ are given by:
\begin{equation}\label{eq:W_mat}
    \dot{q} = \frac12 W(q)\omega = \frac12 \begin{bmatrix} -q_x & -q_y & -q_z \\ q_w & -q_z & q_y \\ q_z & q_w & -q_x \\ -q_y & q_x & q_w \end{bmatrix} \omega
\end{equation}
and the rotation matrix that transforms a vector in the body-fixed satellite frame to the equivalent vector in the global ECI frame is given by \cite{babcock2011cubesat}
\begin{equation}
    \text{rot}(q) = \begin{bmatrix}q_w^{2}+q_x^{2}-q_y^{2}-q_z^{2} & 2\left(q_x q_y-q_w q_z\right) & 2\left(q_x q_z+q_w q_y\right) \\
2\left(q_x q_y+q_w q_z\right) & q_w^{2}-q_x^{2}+q_y^{2}-q_z^{2} & 2\left(q_y q_z-q_w q_x\right) \\
2\left(q_x q_z-q_w q_y\right) & 2\left(q_y q_z+q_w q_x\right) & q_w^{2}-q_x^{2}-q_y^{2}+q_z^{2}  \end{bmatrix}
\end{equation}
\subsection{Satellite Dynamics}
For a magnetorquer-only satellite with an inertia tensor $J$, the angular velocity, given by $\omega$, can be modeled using Euler's equation subject to an input torque $\tau$:
\begin{equation}
    \dot{\omega} = J^{-1}(\tau - \omega \times J\omega)
\end{equation}
It follows that the dynamics of the magnetorquer-only satellite, with $B$ representing the direction local magnetic field, $u$ representing the input dipole of the magnetorquers, and $W(q)$ given in equation \ref{eq:W_mat}, are:
\begin{equation}\label{eq:mtq_dynamics} \dot{x} = f(x, u) = \begin{bmatrix} \dot{\omega} \\ \dot{q} \end{bmatrix} = \begin{bmatrix} J^{-1}(u \times \text{rot}(q)^T B_{ECI} - \omega \times J \omega) \\ \frac12 W(q) \omega \end{bmatrix} \end{equation}
For the 1RW satellite, we augmented the state vector with the reaction wheel momentum $H$, and the dynamics are given by the following equation, where $J_{body}$ is the inertia tensor of the satellite body excluding the reaction wheel, $J_{RW}$ is the inertia tensor of the reaction wheel, $J$ is the inertia tensor of the satellite including the reaction wheel, and $[0\ 1\ 0]^T$ is the axis of rotation of the reaction wheel:
\begin{equation}\label{eq:1RW_dynamics_readable} \dot{x} = f(x, u) = \begin{bmatrix} \dot{\omega} \\ \dot{q} \\ \dot{H} \end{bmatrix} = \begin{bmatrix} J_{body}^{-1}(u \times \text{rot}(q)^T B_{ECI} - \omega \times J \omega) \\ \frac12 W(q) \omega  \\ -u_{RW}-\textrm{diag}(J_{RW})\begin{bmatrix}0\\1\\0\end{bmatrix}\dot{\omega} \end{bmatrix} \end{equation}
Equation \ref{eq:1RW_dynamics_readable} can be rewritten to eliminate $\dot{\omega}$ from the right side of the equation:
\begin{equation}\label{eq:1RW_dynamics} \dot{x} = f(x, u) = \begin{bmatrix} \dot{\omega} \\ \dot{q} \\ \dot{H} \end{bmatrix} = \begin{bmatrix} J_{body}^{-1}(u \times \text{rot}(q)^T B_{ECI} - \omega \times J \omega) \\ \frac12 W(q) \omega  \\ -u_{RW}-\textrm{diag}(J_{RW})\begin{bmatrix}0\\1\\0\end{bmatrix} \left( J_{body}^{-1}(u \times \text{rot}(q)^T B_{ECI} - \omega \times J \omega) \right) \end{bmatrix} \end{equation}
\section{Approach}
The major contribution of this paper is reformulating the trajectory optimization approach proposed for magnetorquer-only control by Gatherer and Manchester \cite{gatherer_magnetorquer-only_2019} to apply to satellites with three magnetorquers and a single reaction wheel by leveraging the dynamics from Equation \ref{eq:1RW_dynamics} and the reaction wheel-specific cost function augmentation described in Equations \ref{eq:angular_momentum_cost} and \ref{eq:stiction_cost}. The combination of a single reaction wheel and trajectory optimization is highly effective and significantly improves upon the pointing accuracy of both magnetorquer-only control with trajectory optimization and attitude control of single reaction wheel systems using traditional approaches like PID control.
\subsection{Trajectory Optimization}
\noindent Given the satellite dynamics given in Equations \ref{eq:mtq_dynamics} and \ref{eq:1RW_dynamics}, we formulated the trajectory optimization problem as follows, with $\ell(x_t, u_t)$ representing the step cost at step $t$, $\ell_f(x_N)$ representing the final cost, and $g(x_t, u_t, B_t)$ representing the forwards integration of the dynamics function $f(x_t, u_t, B_t)$ using the 4th order Runge-Kutta 4 (RK4) method for numerical integration.
\begin{equation} \label{eq:trajopt_formulation}
\begin{aligned}
\minimize_{x_t, u_t} \quad & \ell_f(x_N) + \Sigma_{t=0}^{N-1} \ell(x_t, u_t)\\
\textrm{subject to} \quad & x_{t+1} = g(x_t, u_t, B_t)\\
  &x_0 = x_{\text{init}}    \\
\textrm{with} \quad & t \in \{ 0, 1, ..., N-1 \}
\end{aligned}
\end{equation}
Our cost function is broken down into four subcosts: state cost ($\ell_x$), actuation cost ($\ell_u$), angular momentum cost ($\ell_H$), and stiction cost ($\ell_s$). The state cost is based on the cost function proposed for magnetorquer-only satellites by Gatherer and Manchester\cite{gatherer2018magnetorquer}, with an additional term ($c_{q\omega}$) that provides directional angular velocity damping--increasing cost when the angular velocity is such that the pointing error increases, decreasing cost when the error is decreasing, and zero if it is perpendicular.
\begin{equation}\label{eq:state_cost}
    \ell_x(x_t, q_{goal}, c_\omega, c_q, c_{q\omega}) = \frac12 \omega_t^T \omega_t c_{\omega} + (1 - |q_{goal} \cdot q_t|) c_q - \textrm{sign}(q_{goal} \cdot q_t)q_{goal}W(q_t)\omega c_{q\omega}
\end{equation}

Additionally, we have a control cost, which weights the use of the magnetorquers and the reaction wheel torque. Equation \ref{eq:control_cost} gives this cost -- which is equivalent to $u^T R u$, the quadratic control cost used in linear quadratic regulator (LQR) control, if $R = \frac12 \textrm{diag}([c_{MTQ}\ c_{RW}]^T) c_u$.
\begin{equation}\label{eq:control_cost}
    \ell_u(u_t, c_{MTQ}, c_{RW}, c_u) = \frac12 u_t^T \textrm{diag}([c_{MTQ}\ c_{RW}]^T) u_t c_u
\end{equation}
Finally, there are two additional costs relevant to reaction wheels that are used in the 1RW case. We need to avoid cases where the reaction wheels saturate and where ``stiction'' may occur. When run, the optimization problem has constraints keeping the RW angular momentum beneath its maximum momentum, but we also incur a cost as it gets closer to the angular momentum boundaries--this is to prevent an operational situation where the satellite is trying to fly a trajectory which includes maximized angular momentum in one direction, leaving the satellite unable to correct for disturbances in that direction. However, we do not want to penalize having  angular momentum, so we use a modified quadratic cost. One form of this could be 
\begin{equation}\label{eq:example_angular_momentum_cost}
     \ell_H(x_t, c_{H_1}, c_{H_2}) = \frac12 c_{H_1} \left(\max\left(0,|H|-c_{H_2}\right)\right)^2,
\end{equation}
where $c_{H_1}$ is the weighting of the cost and $c_{H_2}$ is a threshold--when the magnitude of the angular momentum is less than this threshold (when the momentum is either positive or negative), there is no penalty. But once the angular momentum's magnitude passes that threshold (perhaps $\frac23$ of the maximum), there is a quadratic cost. However, the use of the $\max$ function can present issues for the optimization algorithm, as the derivative is not continuous and there is a large flat ``deadzone.'' To address these issues, the ``softplus'' function is used, which approximates the max function continuously and is sometimes used as an activation function in deep neural networks because its derivative is relatively smooth \cite{zheng2015improving}. Thus, adding $c_{H_3}$ as a normalization parameter, the angular momentum cost is:
\begin{equation}\label{eq:angular_momentum_cost}
    \ell_H(x_t, c_{H_1}, c_{H_2}, c_{H_3}) = \frac12 c_{H_1} \left(\frac{\log(1+\exp({(|H|-c_{H_2})c_{H_3}})}{c_{H_3}}\right)^2
\end{equation} 

The other reaction wheel cost is a cost to prevent ``stiction.'' Avoiding saturation in reaction wheels is important, limiting how fast they can spin and how much momentum they can store. But there is also the case where the reaction wheel can spin too slowly. Specifically, if a reaction wheel stops or spins too slowly for too long, there is the risk of it sticking due to friction and internal dynamics. This effect is called ``stiction'' and can prevent a reaction wheel from starting up again. To deal with this (while preventing the issues with non-continuous functions as previously mentioned), we use a small cost based on the ``smoothstep'' function and weighted by $c_{s_1}$ which increases cost if the angular momentum is below a threshold ($c_{s_2}$). Specifically, we use
\begin{equation}\label{eq:stiction_cost}
    \ell_s(x_t, c_{s_1}, c_{s_2}) = \frac12 c_{s_1} \left(\text{smoothstep}\left(\frac{c_{s_2}-|H|}{c_{s_2}} \right)c_{s_2}\right)^2
\end{equation}
where $\text{smoothstep}(x)$ is a smoothed step function:
\begin{equation}
    \text{smoothstep}(x) = \begin{cases}0 \qquad &x \leq 0 \\ 3x^2-2x^2 \qquad &0 < x < 1 \\ 1 \qquad &x \geq 1 \end{cases}
\end{equation}
It should be noted that this cost should be kept low, as it should not be a barrier that prevents the reaction wheel from crossing past zero angular momentum when necessary.
\\ \\
The step cost functions for the magnetorquer-only satellite are:
\begin{align}
    \ell(x_t, u_t, q_{goal}) &= \ell_x(x_t, q_{goal}) + \ell_u(u_t)\\
    \ell_f(x_t, q_{goal}) &= \ell_x(x_t, q_{goal})
\end{align}
and the step costs for the single-reaction-wheel satellite are:
\begin{align}
        \ell(x_t, u_t, q_{goal}) &= \ell_x(x_t, q_{goal}) + \ell_u(u_t) + \ell_H(x_t) + \ell_s(x_t)\\ 
        \ell_f(x_t, q_{goal}) &= \ell_x(x_t, q_{goal}) + \ell_H(x_t) + \ell_s(x_t)
\end{align}
We solve the trajectory optimization problem given in Equation \ref{eq:trajopt_formulation} using our own implementation of ALTRO \cite{howell2019altro}. Among other modifications, we ported ALTRO from Julia to C++, calculated the dynamics and cost Jacobians analytically rather than using autodifferentiation, and eliminated the projection phase of ALTRO to speed up the solver and improve usability on embedded hardware on CubeSats.

\section{Simulation Results}
Using the mathematical formulation outlined above, we evaluated the performance of trajectory optimization for our magnetorquer-only and 1RW systems. The properties of the simulated satellites are given in Table \ref{tab:sat_params}.
\begin{table}[htbp] 
	\fontsize{10}{10}\selectfont
    \caption{Satellite properties (reaction wheel properties below the dashed line).}
        \centering 
    \label{tab:sat_params}
   \begin{tabular}{c | r  } % Column formatting, 
      \hline 
      Property    & Value \\
      \hline 
      Satellite $J_{xx}\ (\text{kg m}^2)$     & 0.005256 \\
      Satellite $J_{yy} = J_{zz}\ (\text{kg m}^2)$& 0.04939  \\
      Maximum magnetic moment (x) (Am$^2$) & 0.19  \\
      \hdashline
      Maximum reaction wheel torque (Nm) & 0.0002 \\
      Momentum storage (Nms) & 0.002 \\
      Reaction wheel moment of inertia (kg m$^2$) & 2e-6 \\
      Reaction wheel rotation axis &  $\ [0\ 1\ 0]^T$ \\
      \hline
   \end{tabular}
\end{table}
\\ \\
As in Manchester \& Gatherer \cite{gatherer2018magnetorquer}, we evaluated the performance of our trajectory planner by planning slews from $q = [0\ 0\ 1\ 0]^T$ to $q = [0\ 1\ 0\ 0]^T$. We also evaluated the performance of our trajectory planner in two other cases: a reduced-attitude goal slew that aims to align the satellite body x-axis with a ground target in the Earth-centered inertial coordinate frame, and a multiple target trajectory that aims to align the satellite body x-axis with a sequence of different ground targets in the Earth-centered inertial coordinate frame.
\\ \\
For all three cases, we initialized the satellite orbital position randomly along a circular ISS orbit (51.5$^\circ$ inclination) with 429 km altitude, and initialized the satellite angular velocity and reaction wheel stored momentum to zero. In the cases for the Reduced Attitude Goals and Multiple Targets, we applied a 20 degree keep-out zone around the sun, to retain realistic constraints that are applicable even in reduced-attitude approaches.

\subsection{Simple 180$^{\circ}$ Slew}
\begin{figure}[htp]
	\centering
	\begin{subfigure}
	\centering
	\includegraphics[width=0.45\linewidth]{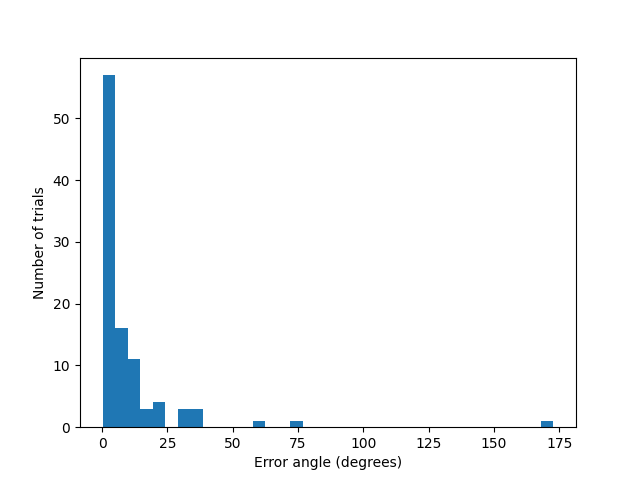}
	\end{subfigure}
	\begin{subfigure}
	\centering
	\includegraphics[width=0.45\textwidth]{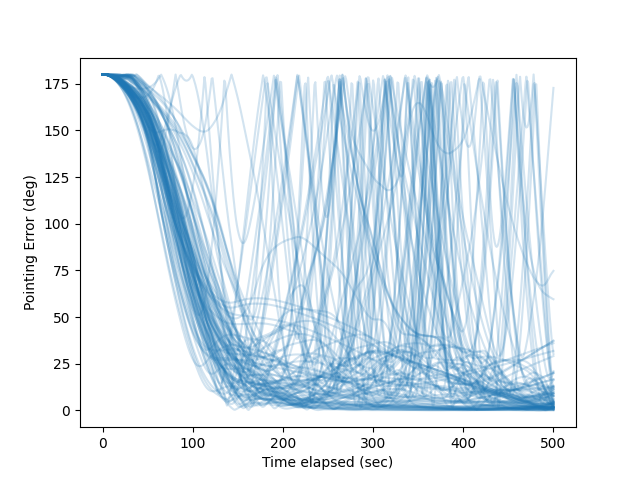}
  \end{subfigure}
  \caption{Results of 100-trial Monte Carlo simulation for planning 180$^{\circ}$ slews for a magnetorquer-only satellite. Left: Final attitude error for all trials. Right: Visualization of attitude error over time for all trials.}
  \label{fig:mtq_montecarlo}
\end{figure}

\begin{figure}[htp]
	\centering
	\begin{subfigure}
	\centering
	\includegraphics[width=0.45\linewidth]{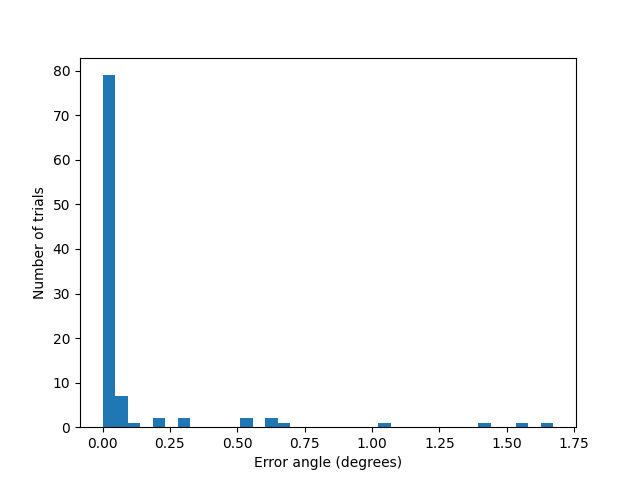}
	\end{subfigure}
	\begin{subfigure}
	\centering
	\includegraphics[width=0.45\textwidth]{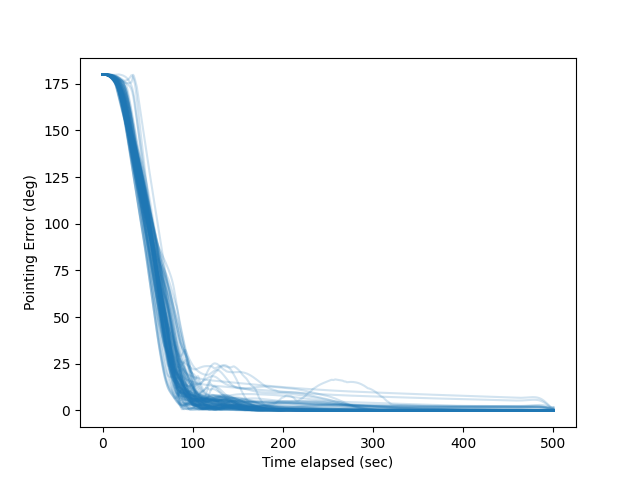}
  \end{subfigure}
  \caption{Results of 100-trial Monte Carlo simulation for planning 180$^{\circ}$ slews for a 1RW satellite. Left: Final attitude error for all trials. Right: Visualization of attitude error over time for all trials.}
  \label{fig:1RW_montecarlo}
\end{figure}

\begin{figure}[htp]
	\centering
	\begin{subfigure}
	\centering
	\includegraphics[width=0.48\textwidth]{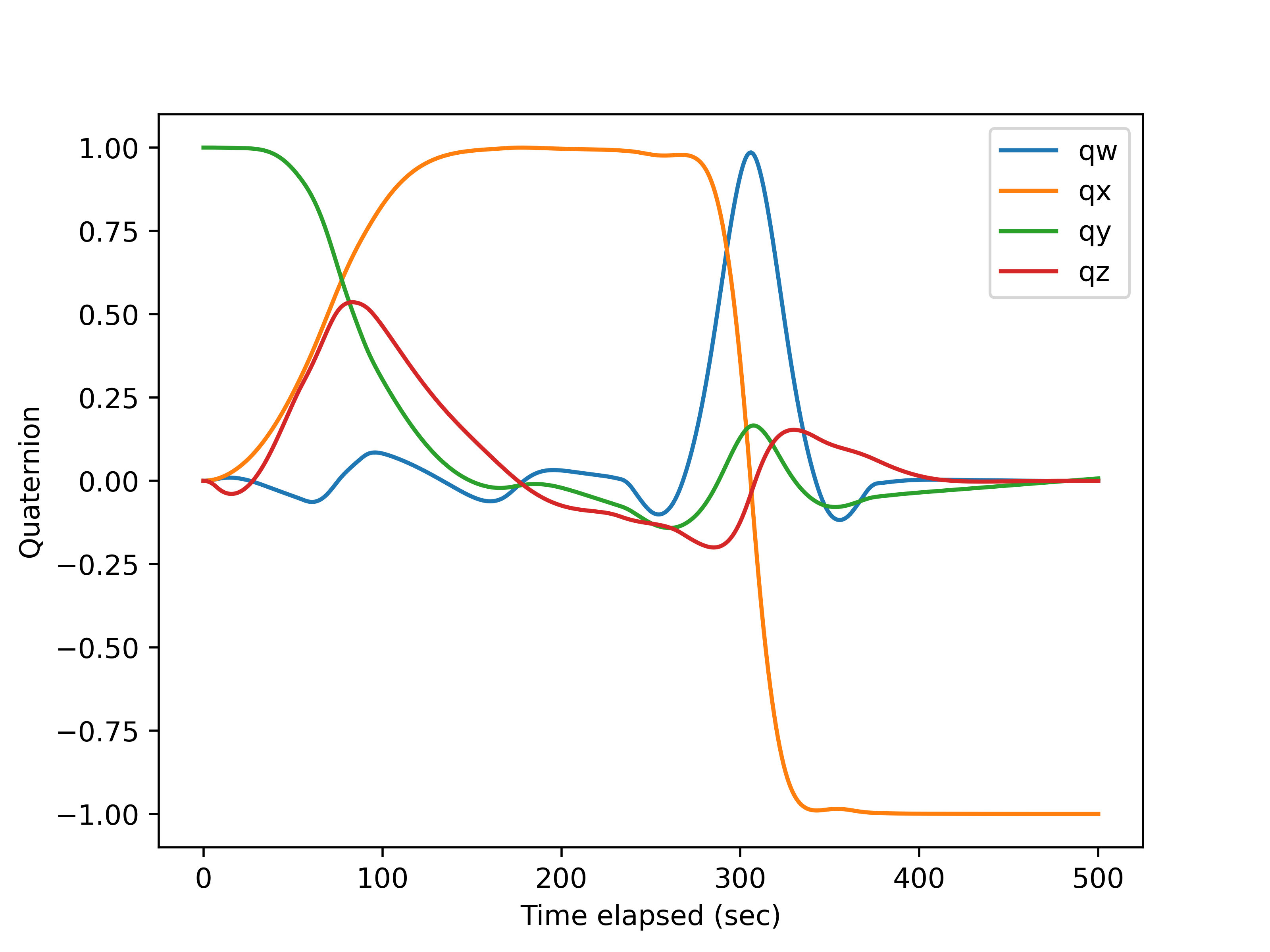}
	\end{subfigure}
	\begin{subfigure}
	\centering
	\includegraphics[width=0.48\textwidth]{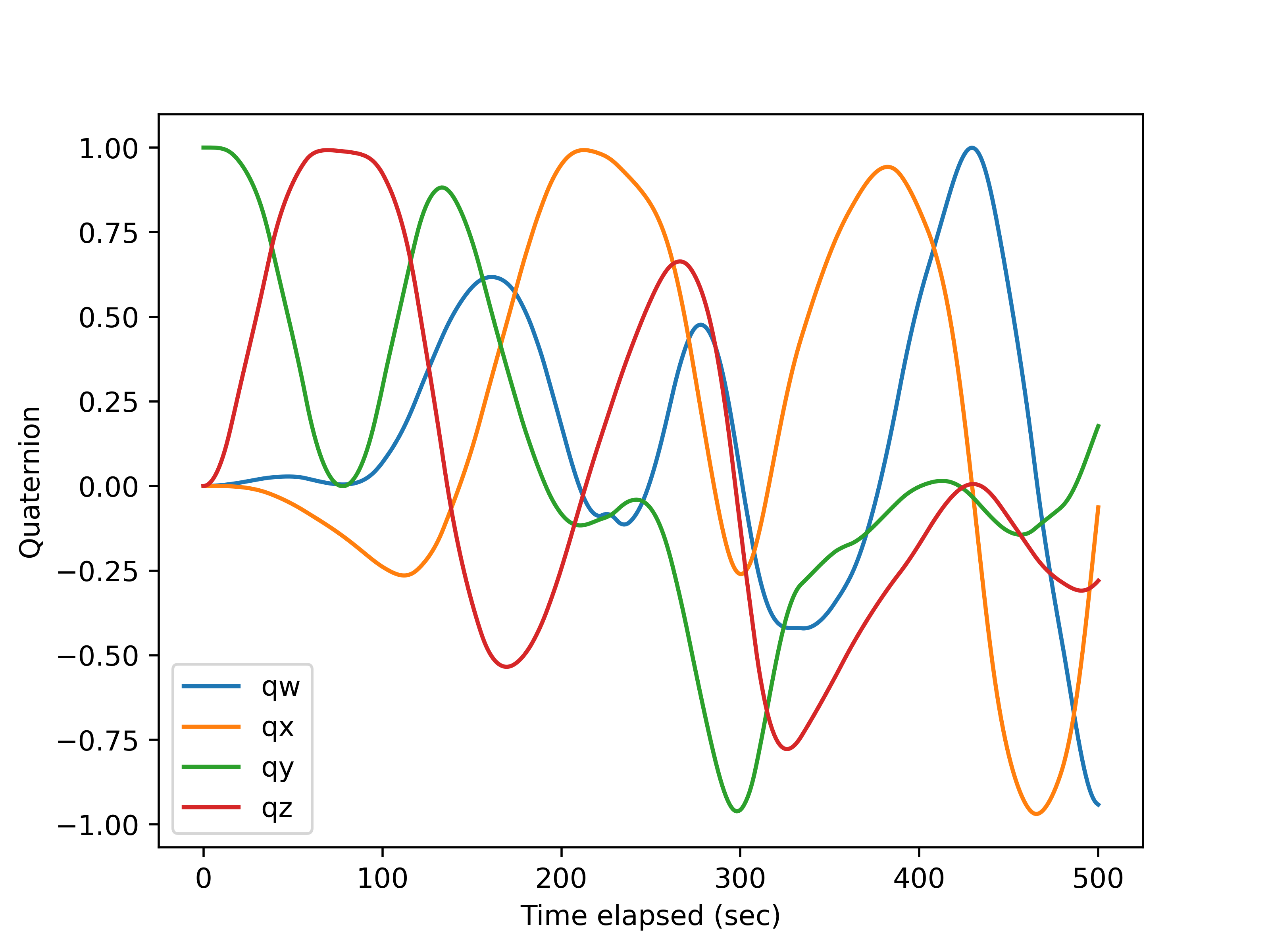}
  \end{subfigure}
  \caption{Two planned slews for a magnetorquer-only satellite. Left: ``Good convergence" case, demonstrating a slew to target. Right: ``Bad convergence" case, showing noisy oscillation in attitude and no convergence to the target.}
  \label{fig:mtq_quat_trajectories}
\end{figure}

\begin{figure}[htp]
	\centering
	\begin{subfigure}
	\centering
	\includegraphics[width=0.45\textwidth]{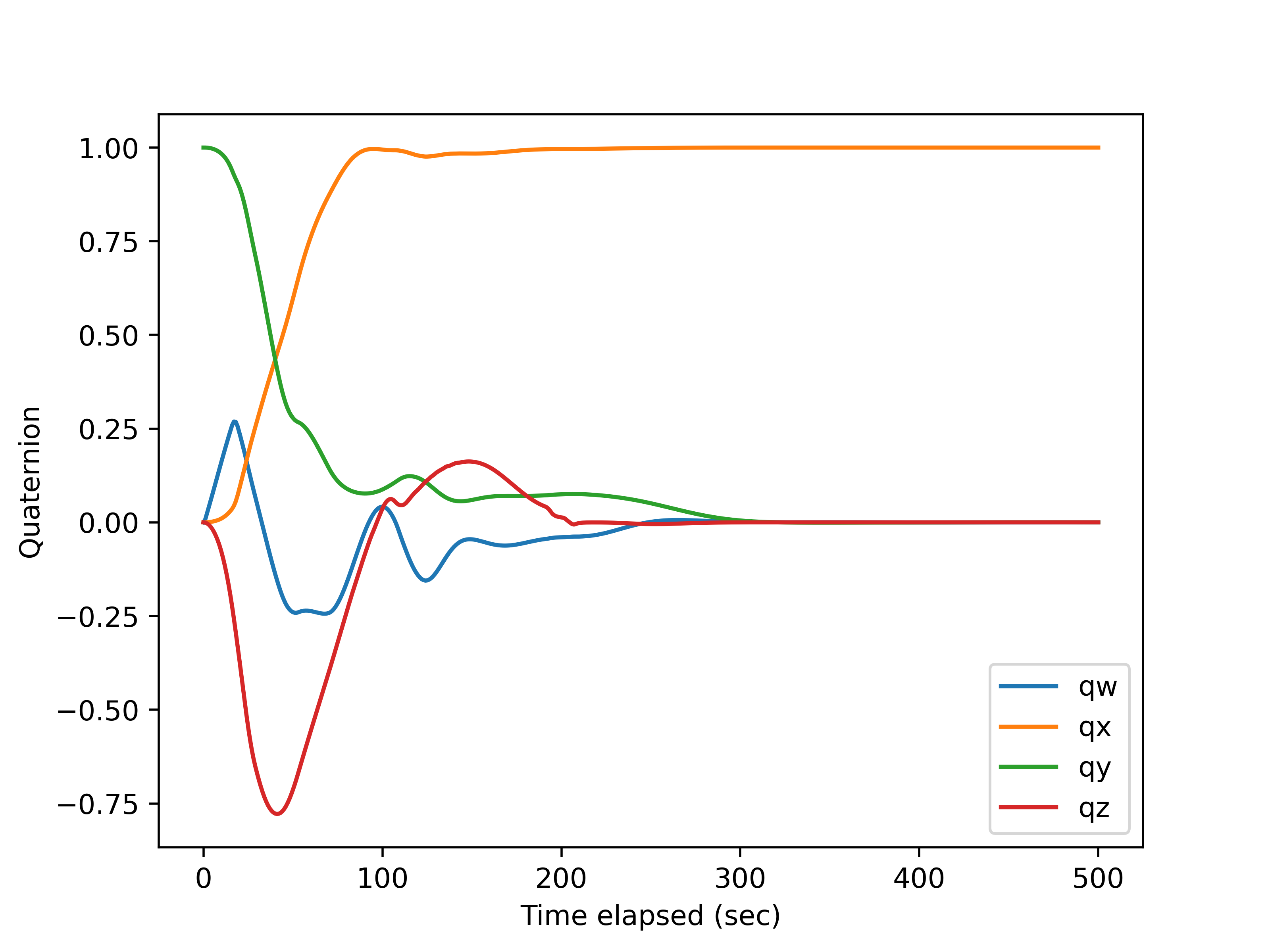}
	\end{subfigure}
	\begin{subfigure}
	\centering
	\includegraphics[width=0.45\textwidth]{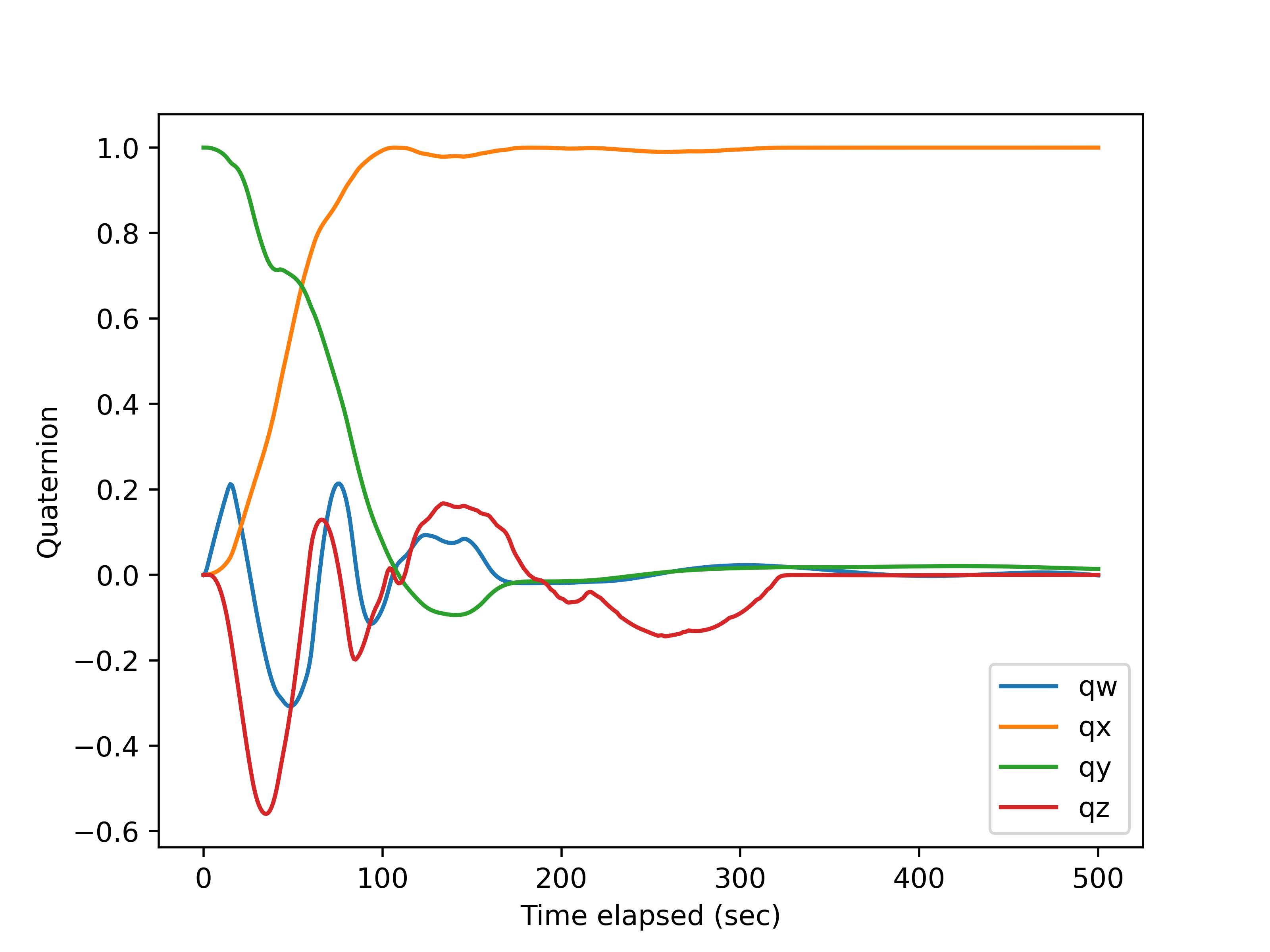}
  \end{subfigure}
  \caption{Two planned slews for a 1RW satellite. Left: ``Good convergence" case, demonstrating a direct slew to target. Right: ``Bad convergence" case, demonstrating a slow oscillation towards the target.}
  \label{fig:1RW_quat_trajectories}
\end{figure}

\begin{figure}[htp]
	\centering
	\includegraphics[width=0.45\textwidth]{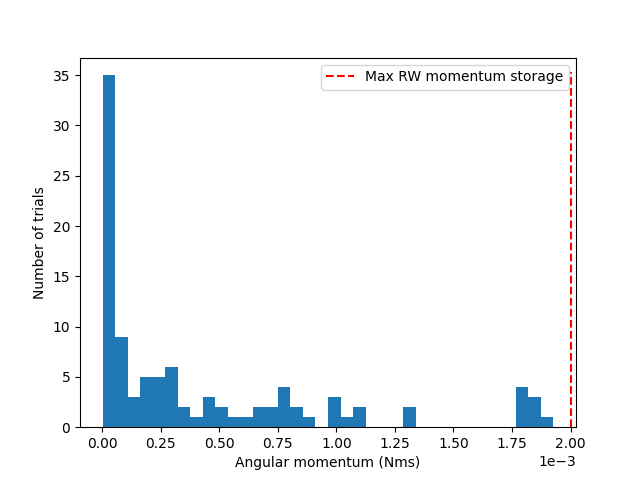}
    \caption{Results of 100-trial Monte Carlo simulation for planning 180$^{\circ}$ slews for a 1RW satellite: final RW stored angular momentum for all trials.}
  \label{fig:1RW_momentum}
\end{figure}

We performed a 100-trial Monte Carlo simulation over the randomized initial orbital conditions in order to analyze the final attitude error -- the angle between the final attitude quaternion from the closed-loop simulation and $[0\ 1\ 0\ 0]^T$.
\\ \\
 We show the final attitude error for 100 trials of planning and executing $180^{\circ}$ slews for a magnetorquer-only satellite and visualize the final attitude error over time for all 100 trials in Figure \ref{fig:mtq_montecarlo}. In some cases, the trajectory planning algorithm converges to a poor trajectory, resulting in rapid spins and final attitude error on the order of tens of degrees. Overall, the final attitude error has a mean of 10.21 degrees, a median of 4.20 degrees, and a standard deviation of 20.21 degrees.
\\ \\
For the 1RW satellite,  we show the final attitude error for 100 trials of planning and executing $180^{\circ}$ slews and visualize attitude error over time for all 100 trials in Figure \ref{fig:1RW_montecarlo}. As with the magnetorquer-only satellite, sometimes the trajectory planning algorithm converges to a suboptimal trajectory, but for the 1RW satellite, this results in final pointing errors between 1 and 10 degrees. The 1RW slews have a lower mean (0.11 degrees) and median (1.2e-4 degrees) final attitude error than the magnetorquer-only slews, and, crucially, have a lower variance. The standard deviation of final attitude error for the 1RW slews is just 0.31 degrees.
\\ \\
In Figure \ref{fig:mtq_quat_trajectories} we show the planned quaternion attitude over time for a ``good convergence" case and for a ``bad convergence" case with high final attitude error, both for a magnetorquer-only satellite. The planned attitude in the ``bad" convergence case shown in Figure \ref{fig:mtq_quat_trajectories} shows oscillation that never approaches the target. This drift is characteristic of the planned magnetorquer-only trajectories that result in more than $10^{\circ}$ degrees of final attitude error.
\\ \\
For the 1RW satellite, the planned quaternion attitude over time for a ``good" and ``bad" convergence case is shown in Figure \ref{fig:1RW_quat_trajectories}. The ``bad convergence" case shown here, like the one shown in Figure \ref{fig:mtq_quat_trajectories}, slowly approaches the target, but in this case oscillates towards the goal, ultimately achieving less than $2^{\circ}$ of final attitude error, which would still be acceptable for many CubeSat missions. In all 100 trials, the 1RW satellite had less than $2^\circ$ of final attitude error.
\\ \\
We show the final angular momentum of the reaction wheel for 100 trials of planning and executing $180^{\circ}$ slews for a 1RW satellite in Figure \ref{fig:1RW_momentum}. Notably, the maximum final angular momentum is 1.9e-3 Nms, or $95\%$ of the reaction wheel momentum storage, and as such, the reaction wheel is never fully saturated at the end of the trajectory. The mean final angular momentum is 4.2e-4 Nms, the median final angular momentum is 2.0e-5 Nms, and the standard deviation of the final angular momentum is 5.3e-4 Nms.
\subsection{Reduced Attitude Goal Slews}
\begin{figure}[htp]
	\centering
	\begin{subfigure}
	\centering
	\includegraphics[width=0.45\linewidth]{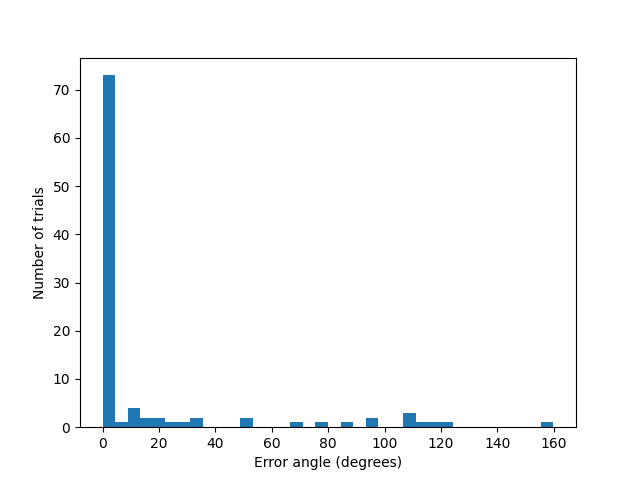}
	\end{subfigure}
	\begin{subfigure}
	\centering
	\includegraphics[width=0.45\textwidth]{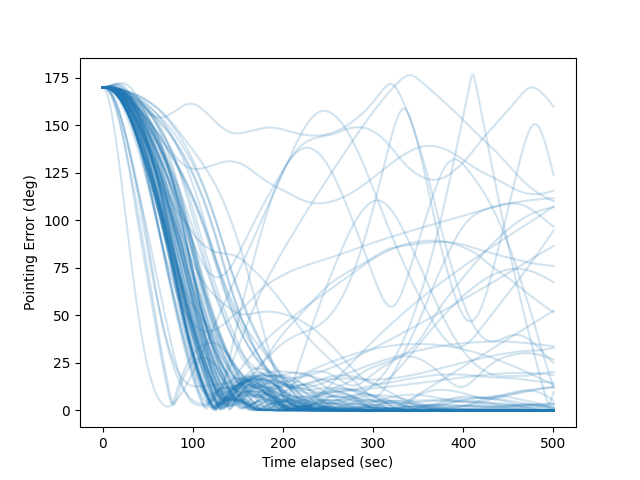}
  \end{subfigure}
  \caption{Results of 100-trial Monte Carlo simulation for planning reduced attitude goal slews for a magnetorquer-only satellite. Left: Final attitude error for all trials. Right: Visualization of attitude error over time for all trials.}
  \label{fig:mtq_montecarlo_imaging}
\end{figure}

\begin{figure}[htp]
	\centering
	\begin{subfigure}
	\centering
	\includegraphics[width=0.45\linewidth]{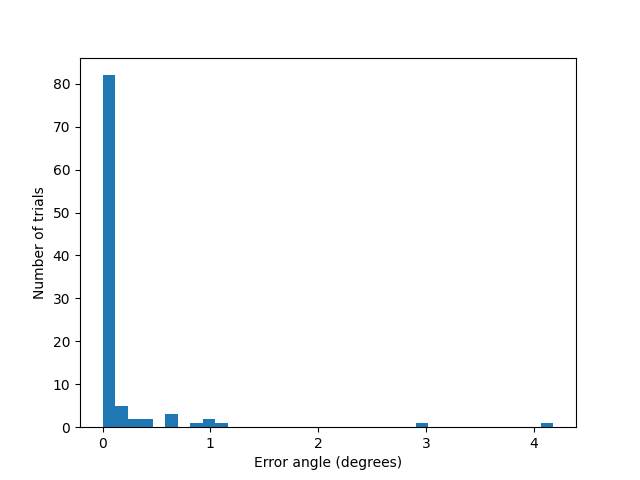}
	\end{subfigure}
	\begin{subfigure}
	\centering
	\includegraphics[width=0.45\textwidth]{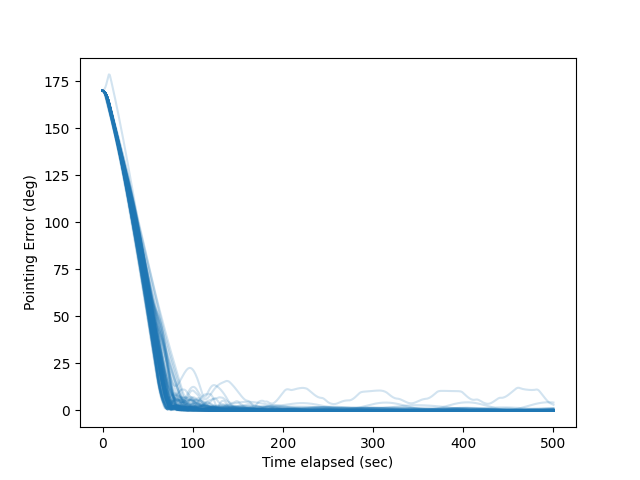}
  \end{subfigure}
  \caption{Results of 100-trial Monte Carlo simulation for planning reduced attitude goal slews for a 1RW satellite. Left: Final attitude error for all trials. Right: Visualization of attitude error over time for all trials.}
  \label{fig:1RW_montecarlo_imaging}
\end{figure}

\begin{figure}[htp]
	\centering
	\begin{subfigure}
	\centering
	\includegraphics[width=0.48\textwidth]{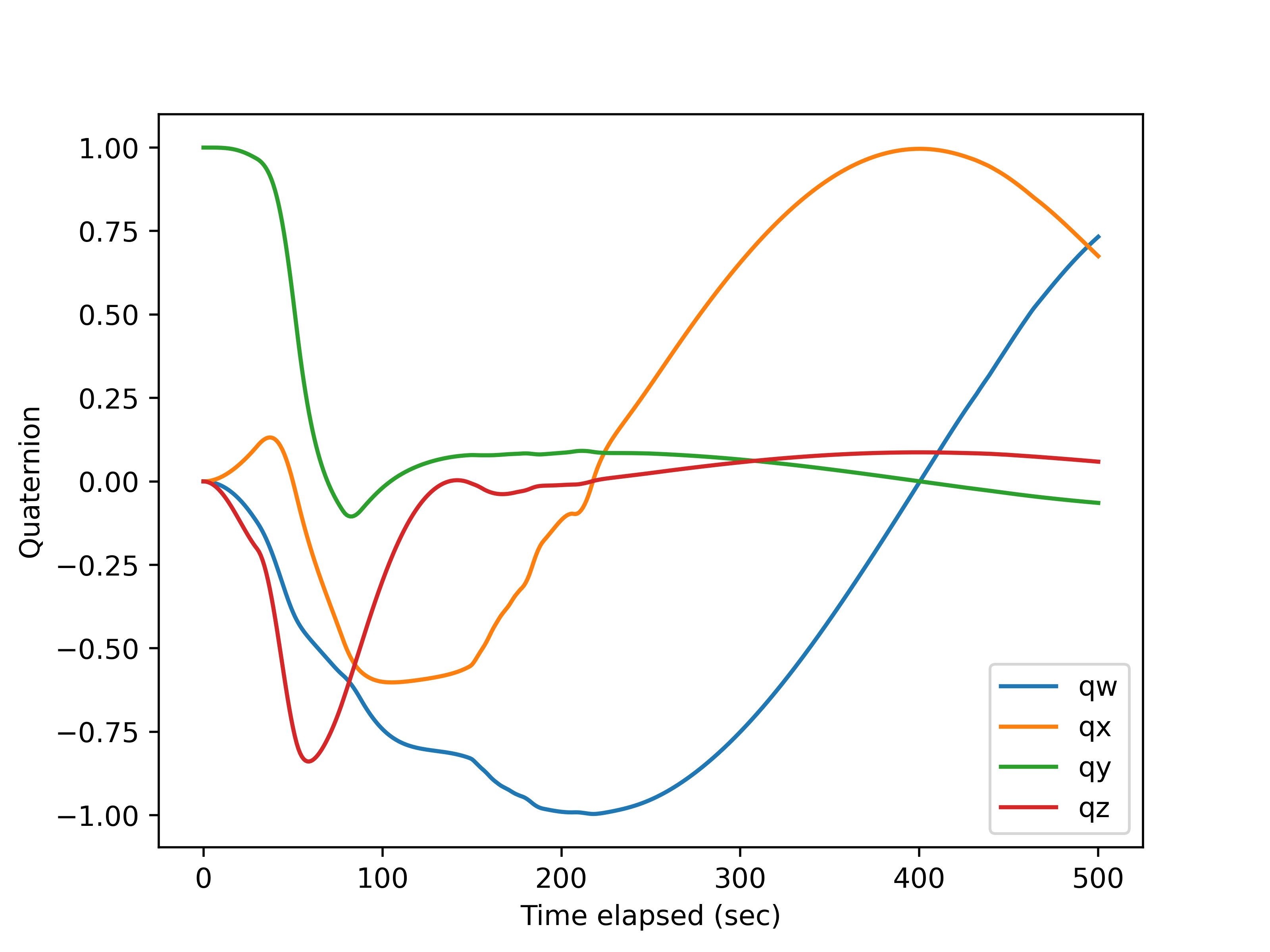}
	\end{subfigure}
	\begin{subfigure}
	\centering
	\includegraphics[width=0.48\textwidth]{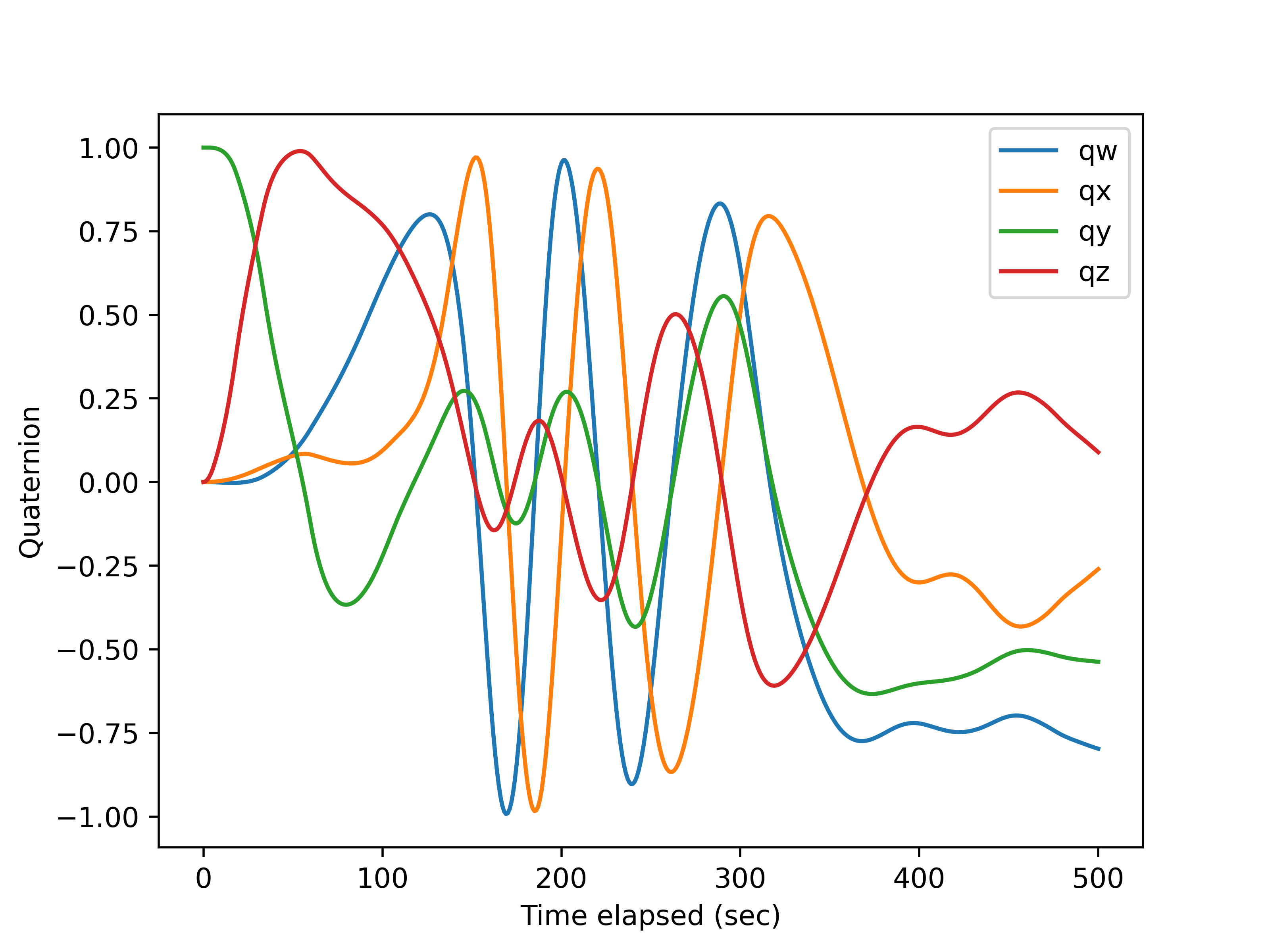}
  \end{subfigure}
  \caption{Two planned slews for a magnetorquer-only satellite. Left: ``Good convergence" case, demonstrating a slew to target. Right: ``Bad convergence" case, showing noisy oscillation in attitude and no convergence to the target.}
  \label{fig:mtq_quat_trajectories_imaging}
\end{figure}

\begin{figure}[htp]
	\centering
	\begin{subfigure}
	\centering
	\includegraphics[width=0.45\textwidth]{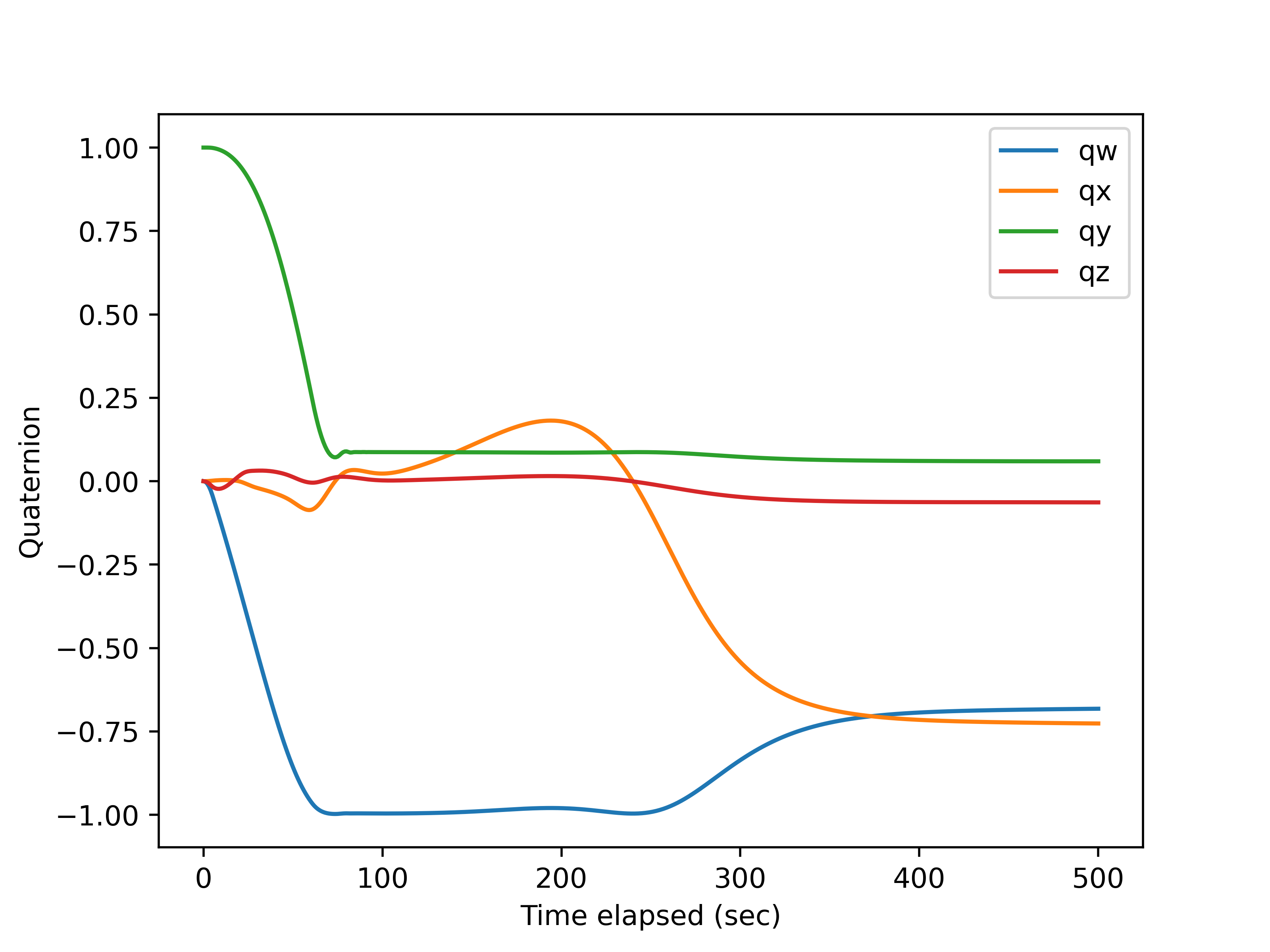}
	\end{subfigure}
	\begin{subfigure}
	\centering
	\includegraphics[width=0.45\textwidth]{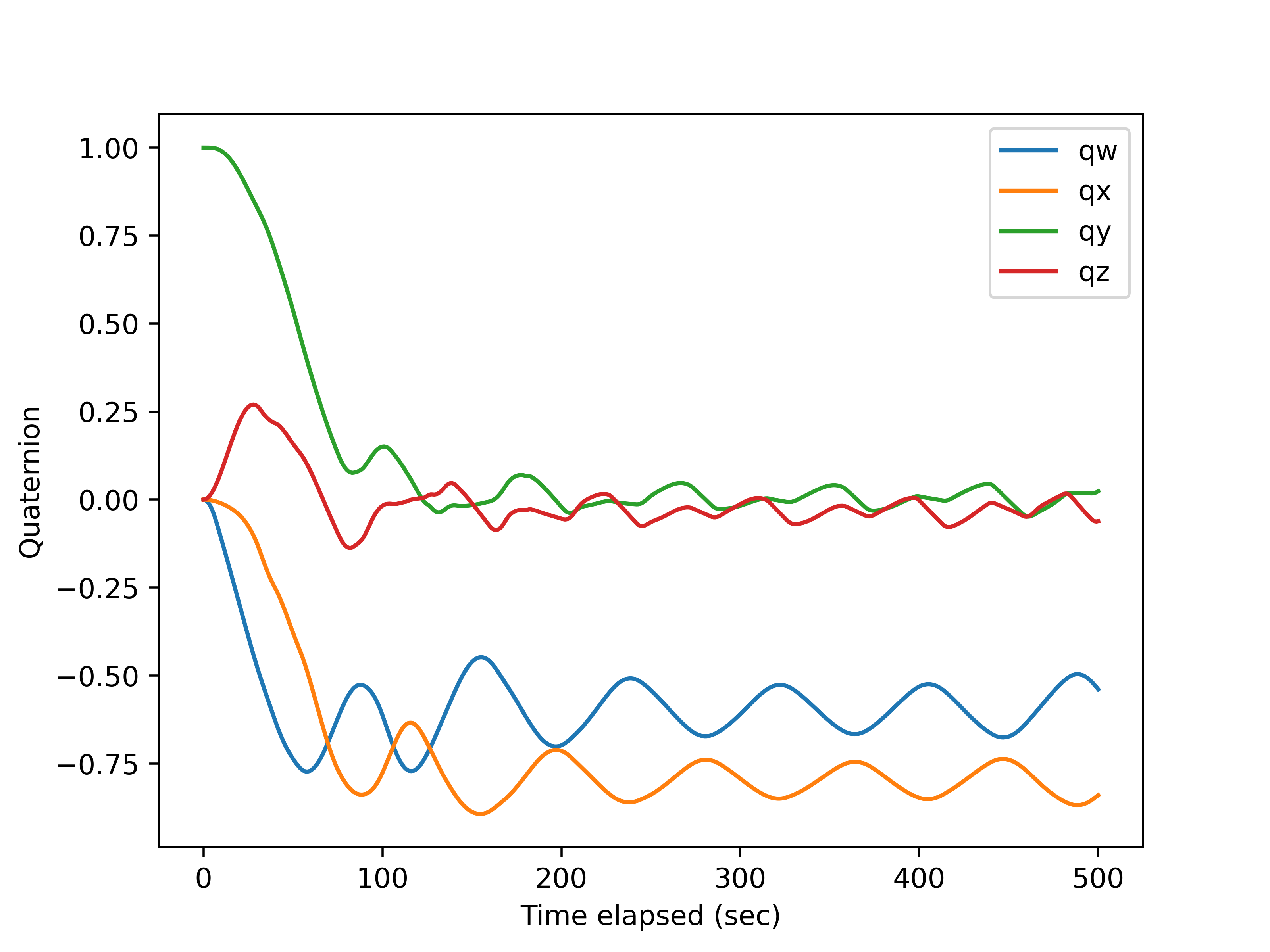}
  \end{subfigure}
  \caption{Two planned slews for a 1RW satellite. Left: ``Good convergence" case, demonstrating a direct slew to target. Right: ``Bad convergence" case, demonstrating a slow oscillation around the target.}
  \label{fig:1RW_quat_trajectories_imaging}
\end{figure}

\begin{figure}[htp]
	\centering
	\includegraphics[width=0.45\textwidth]{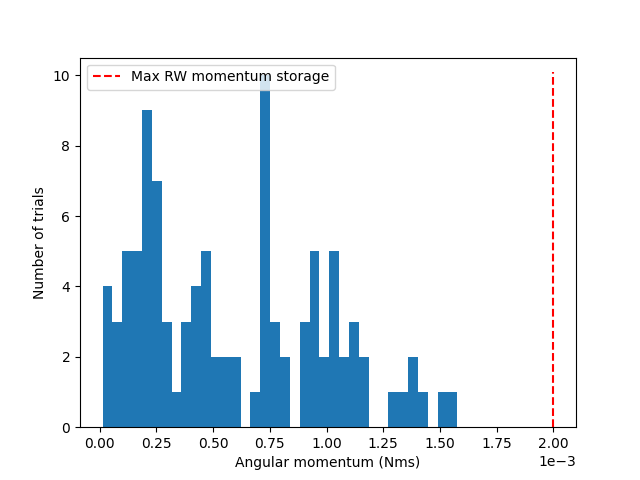}
    \caption{Results of 100-trial Monte Carlo simulation for planning reduced attitude goal slews for a 1RW satellite: final RW stored angular momentum for all trials.}
  \label{fig:1RW_momentum_imaging}
\end{figure}

For the reduced attitude goal slew, we planned slews to align the satellite body x-axis with $[\cos(10^\circ), 0, \sin(10^\circ)]^T$ in the Earth-centered inertial (ECI) frame. Notably, the reduced attitude goals have a simplified aim: to align a specific satellite body axis with a target vector in the ECI coordinate frame. As such, they do not specify a full three-dimensional target orientation. Given an orientation with body x-axis aligned with the target, any rotation around the body x-axis yields another valid orientation. (The equations describing the set of quaternions which meet this goal can be succinctly described.\cite{McKeen_Method_2022})  In comparison, the simple $180^\circ$ slews specify a single valid target-pointing orientation. Due to the larger set of possibilities, the magnetorquer-only and 1RW trajectories should both achieve lower pointing error on the reduced attitude goal slew than the simple $180^\circ$ slew.
\\ \\
As for the simple $180^\circ$ slews, we performed a 100-trial Monte Carlo simulation over the randomized initial orbital conditions in order to analyze the final attitude error -- the angle between the satellite body x-axis and $[\cos(10^\circ), 0, \sin(10^\circ)]^T$ in the ECI frame at the end of the trajectory.
\\ \\
 We show the final attitude error for 100 trials of planning and executing reduced attitude goal slews for a magnetorquer-only satellite and visualize the final attitude error over time for all 100 trials in Figure \ref{fig:mtq_montecarlo_imaging}. In some cases, the trajectory planning algorithm converges to a poor trajectory, resulting in rapid spins and final attitude error on the order of tens of degrees. Overall, the final attitude error has a mean of 16.16 degrees, a median of 0.04 degrees, and a standard deviation of 34.99 degrees.
\\ \\
For the 1RW satellite,  we show the final attitude error for 100 trials of planning and executing reduced attitude goal slews and visualize attitude error over time for all 100 trials in Figure \ref{fig:1RW_montecarlo_imaging}. As with the magnetorquer-only satellite, sometimes the trajectory planning algorithm converges to a suboptimal trajectory, but for the 1RW satellite, this results in final pointing errors between 1 and 10 degrees. The 1RW slews have a lower mean (0.16 degrees) and median (2.6e-3 degrees) final attitude error than the magnetorquer-only slews, and, crucially, have a lower variance. The standard deviation of final attitude error for the 1RW slews is just 0.55 degrees.
\\ \\
In Figure \ref{fig:mtq_quat_trajectories_imaging} we show the planned quaternion attitude over time for a ``good convergence" case and for a ``bad convergence" case with high final attitude error, both for a magnetorquer-only satellite. The planned attitude in the ``bad" convergence case shown in Figure \ref{fig:mtq_quat_trajectories_imaging} shows oscillation that slows but never approaches the target. This drift is characteristic of the planned magnetorquer-only trajectories that result in more than $10^{\circ}$ degrees of final attitude error.
\\ \\
For the 1RW satellite, the planned quaternion attitude over time for a ``good" and ``bad" convergence case is shown in Figure \ref{fig:1RW_quat_trajectories_imaging}. The ``bad convergence" case shown here, like the one shown in Figure \ref{fig:mtq_quat_trajectories_imaging}, oscillates, but in this case oscillates around the goal, ultimately achieving less than $5^{\circ}$ of final attitude error, which would still be acceptable for many CubeSat missions. In all 100 trials, the 1RW satellite had less than $5^\circ$ of final attitude error.
\\ \\
We show the final angular momentum of the reaction wheel for 100 trials of planning and executing $180^{\circ}$ slews for a 1RW satellite in Figure \ref{fig:1RW_momentum_imaging}. Notably, the maximum final angular momentum is 1.6e-3 Nms, or $80\%$ of the reaction wheel momentum storage, and as such, the reaction wheel is never fully saturated at the end of the trajectory. The mean final angular momentum is 5.9e-4 Nms, the median final angular momentum is 5.0e-4 Nms, and the standard deviation of the final angular momentum is 4.1e-4 Nms.

\subsection{Multiple Targets Slew}
\begin{figure}[htp]
	\centering
	\begin{subfigure}
	\centering
	\includegraphics[width=0.45\linewidth]{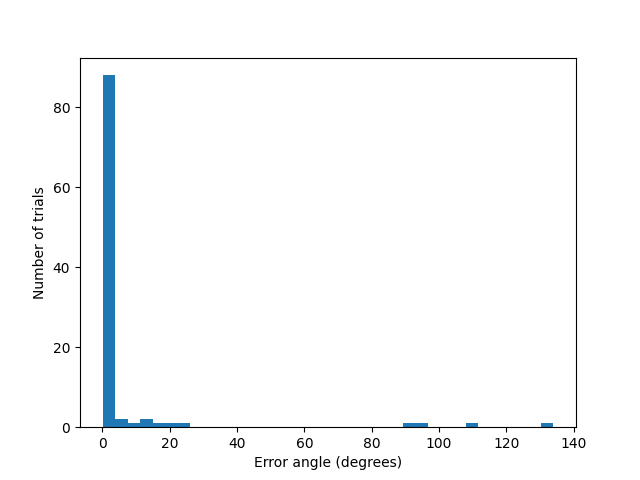}
	\end{subfigure}
	\begin{subfigure}
	\centering
	\includegraphics[width=0.45\textwidth]{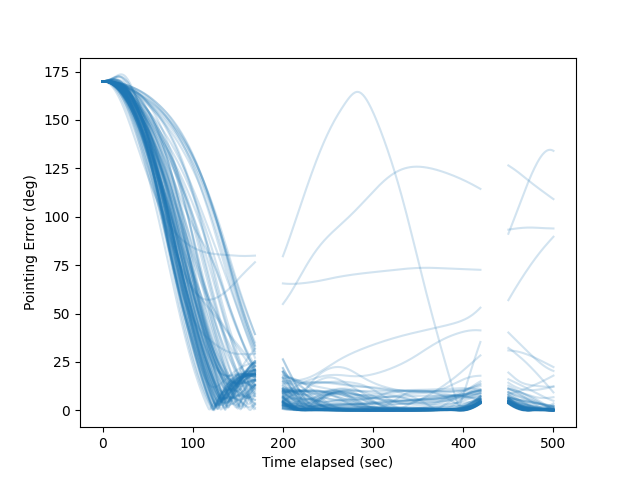}
  \end{subfigure}
  \caption{Results of 100-trial Monte Carlo simulation for planning multiple target trajectories for a magnetorquer-only satellite. Left: Final attitude error for all trials. Right: Visualization of attitude error over time for all trials.}
  \label{fig:mtq_montecarlo_multi}
\end{figure}

\begin{figure}[htp]
	\centering
	\begin{subfigure}
	\centering
	\includegraphics[width=0.45\linewidth]{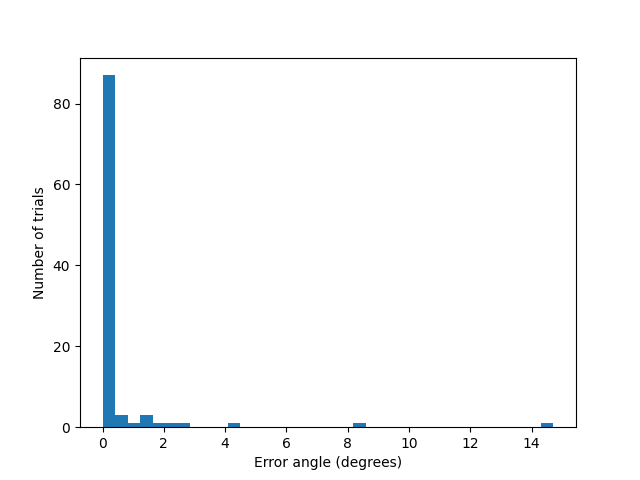}
	\end{subfigure}
	\begin{subfigure}
	\centering
	\includegraphics[width=0.45\textwidth]{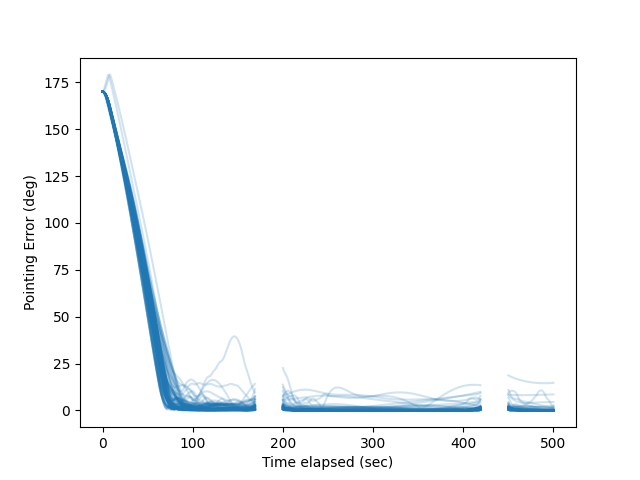}
  \end{subfigure}
  \caption{Results of 100-trial Monte Carlo simulation for planning multiple target trajectories for a 1RW satellite. Left: Final attitude error for all trials. Right: Visualization of attitude error over time for all trials.}
  \label{fig:1RW_montecarlo_multi}
\end{figure}

\begin{figure}[htp]
	\centering
	\begin{subfigure}
	\centering
	\includegraphics[width=0.48\textwidth]{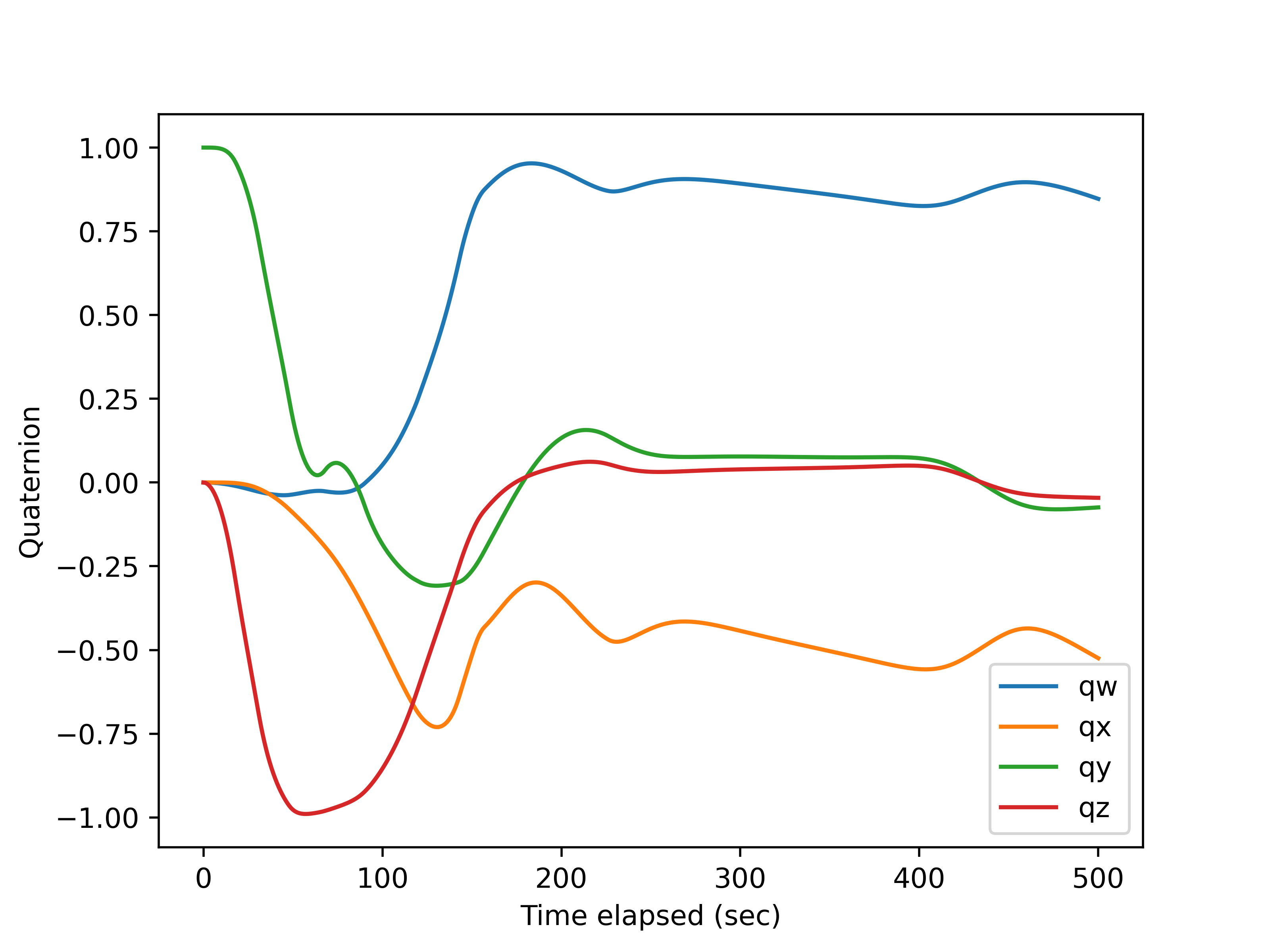}
	\end{subfigure}
	\begin{subfigure}
	\centering
	\includegraphics[width=0.48\textwidth]{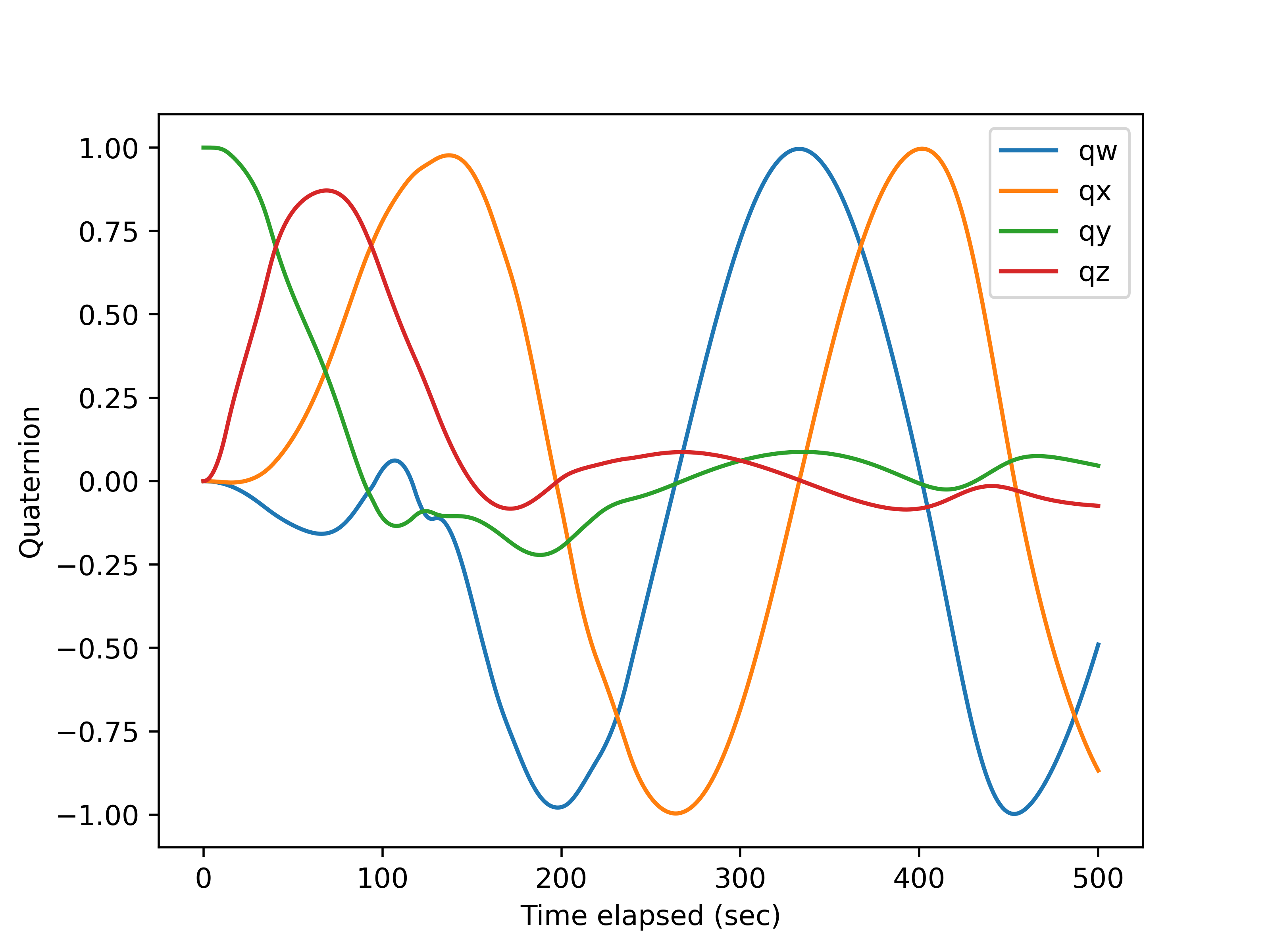}
  \end{subfigure}
  \caption{Two planned slews for a magnetorquer-only satellite. Left: ``Good convergence" case, demonstrating a slew to multiple targets. Right: ``Bad convergence" case, showing noisy oscillation in attitude and no convergence to the target.}
  \label{fig:mtq_quat_trajectories_multi}
\end{figure}

\begin{figure}[htp]
	\centering
	\begin{subfigure}
	\centering
	\includegraphics[width=0.45\textwidth]{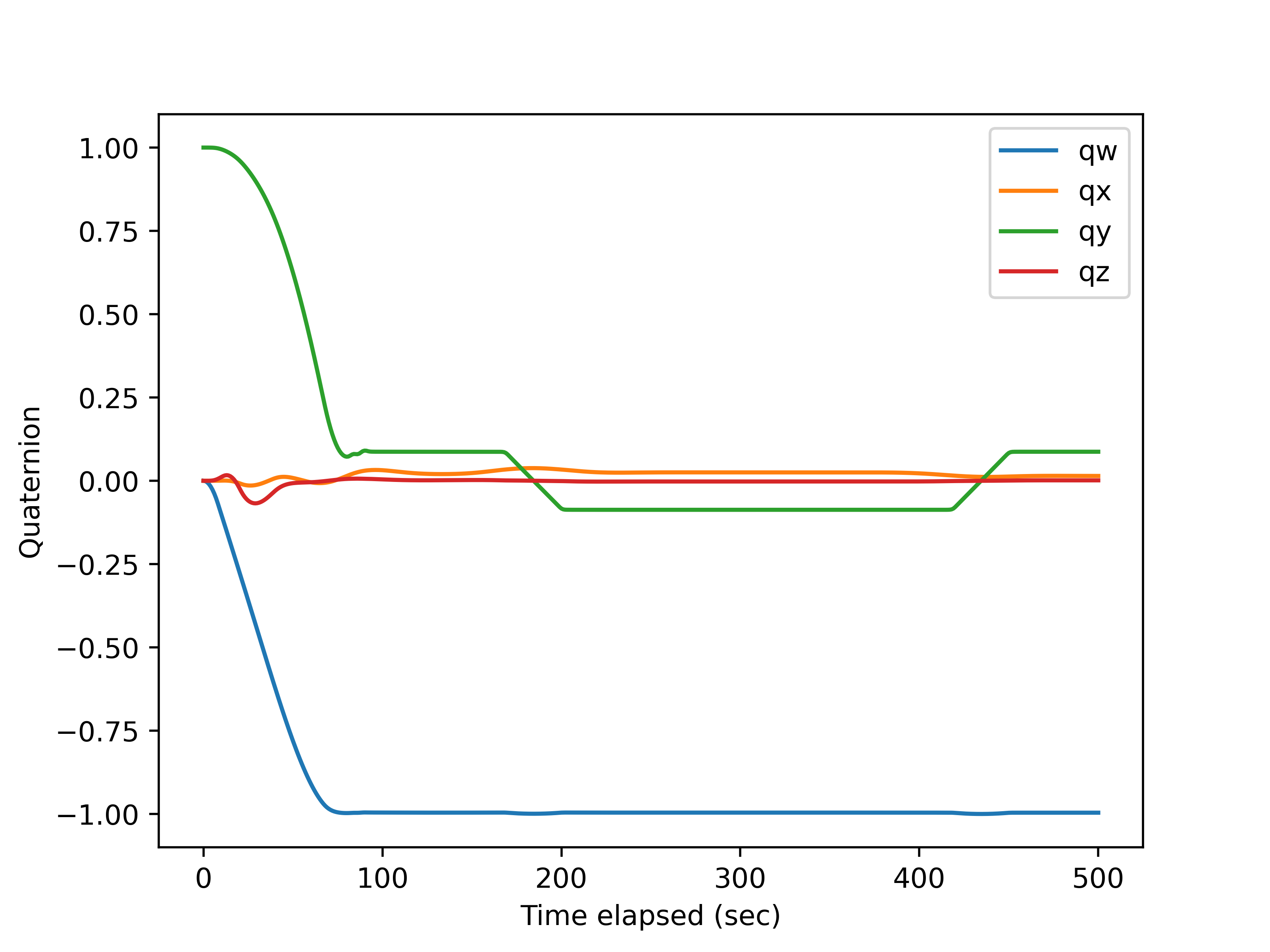}
	\end{subfigure}
	\begin{subfigure}
	\centering
	\includegraphics[width=0.45\textwidth]{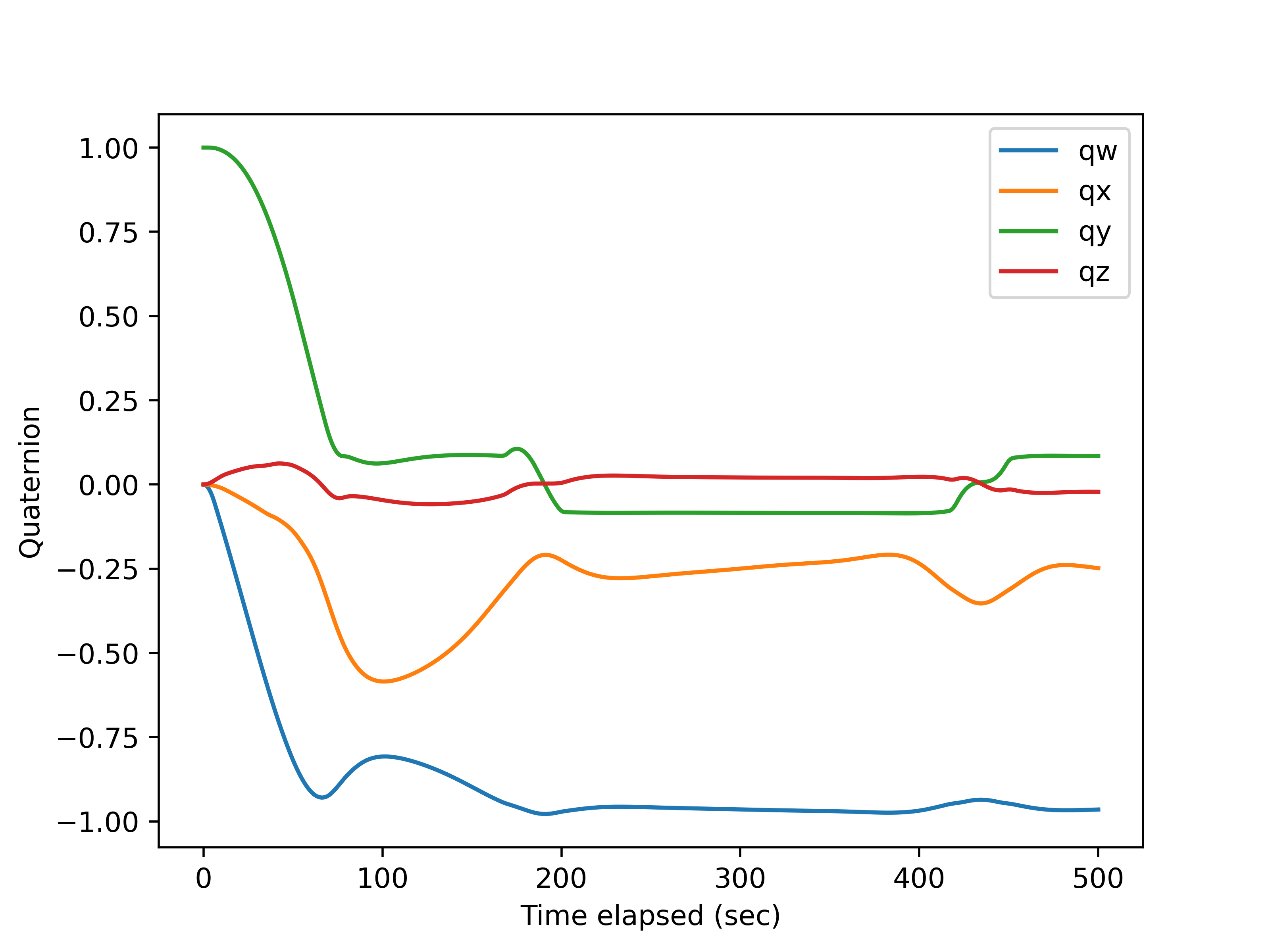}
  \end{subfigure}
  \caption{Two planned slews for a 1RW satellite. Left: ``Good convergence" case, demonstrating a slew to multiple targets. Right: ``Bad convergence" case, demonstrating only approximate convergence to the target.}
  \label{fig:1RW_quat_trajectories_multi}
\end{figure}

For the multiple target trajectory, we planned slews to align the satellite body x-axis with $[\cos(10^\circ), 0, \sin(10^\circ)]^T$ in the Earth-centered inertial (ECI) frame for the first 170 seconds of the slew, with $[\cos(-10^\circ), 0, \sin(-10^\circ)]^T$ in the ECI frame from 200 to 420 seconds, and with $[\cos(10^\circ), 0, \sin(10^\circ)]^T$ in the ECI frame from 450 to 500 seconds. During the multiple target trajectories, the cost function does not penalize deviation from the target quaternion from 170 to 200 seconds and from 420 to 450 seconds -- these periods are intended to be for transitioning between targets, not for target tracking. During this ``free time,'' the planner does not consider the pointing score, so that the period can be used for motions that transition to the next goal. Additionally, like the reduced attitude goal slews, the multiple target trajectories seek only to align a single satellite body axis with an ECI target, and do not fully specify a required satellite orientation.
As for the simple $180^\circ$ slews, we performed a 100-trial Monte Carlo simulation over the randomized initial orbital conditions in order to analyze the final attitude error -- the angle between the satellite body x-axis and $[\cos(10^\circ), 0, \sin(10^\circ)]^T$ in the ECI frame at the end of the trajectory.
\\ \\
 We show the final attitude error for 100 trials of planning and executing multiple target trajectories for a magnetorquer-only satellite and visualize the final attitude error over time for all 100 trials in Figure \ref{fig:mtq_montecarlo_multi}. In some cases, the trajectory planning algorithm converges to a poor trajectory, resulting in rapid spins and final attitude error on the order of tens of degrees. Overall, the final attitude error has a mean of 5.70 degrees, a median of 0.22 degrees, and a standard deviation of 21.27 degrees. The attitude error at the end of the period pointing at the second target has a mean of 9.02 degrees, a median of 4.77 degrees, and a standard deviation of 14.45 degrees. The attitude error at the end of the period pointing at the first target has a mean of 19.40 degrees, a median of 18.59 degrees, and a standard deviation of 11.38 degrees.
\\ \\
For the 1RW satellite,  we show the final attitude error for 100 trials of planning and executing reduced attitude goal slews and visualize attitude error over time for all 100 trials in Figure \ref{fig:1RW_montecarlo_multi}. As with the magnetorquer-only satellite, sometimes the trajectory planning algorithm converges to a suboptimal trajectory, but for the 1RW satellite, this usually results in final pointing errors between 1 and 10 degrees. The 1RW slews have a lower mean (0.45 degrees) and median (0.03 degrees) final attitude error than the magnetorquer-only slews, and, crucially, have a lower variance. The standard deviation of final attitude error for the 1RW slews is just 1.77 degrees. At the end of the period pointing at the second target, the pointing error has a mean of 1.53 degrees, a median of 0.78 degrees, and a standard deviation of 2.11 degrees. At the end of the period pointing at the first target, the pointing error has a mean of 2.22 degrees, a median of 1.60 degrees, and a standard deviation of 2.37 degrees.
\\ \\
In Figure \ref{fig:mtq_quat_trajectories_multi} we show the planned quaternion attitude over time for a ``good convergence" case and for a ``bad convergence" case with high final attitude error, both for a magnetorquer-only satellite. The planned attitude in the ``bad" convergence case shown in Figure \ref{fig:mtq_quat_trajectories_multi} shows oscillation that slows but never approaches the target. This oscillation is characteristic of the planned magnetorquer-only trajectories that result in more than $10^{\circ}$ degrees of final attitude error. Even in the ``good convergence" case, the magnetorquer-only satellite only achieves approximate convergence to each target, and is not stationary during any of the three imaging periods.
\\ \\
For the 1RW satellite, the planned quaternion attitude over time for a ``good" and ``bad" convergence case is shown in Figure \ref{fig:1RW_quat_trajectories_multi}. The ``bad convergence" case shown here roughly converges to each target, but still is not stationary during any of the three imaging periods. In this case, the 1RW satellite ultimately achieves less than $15^{\circ}$ of final attitude error, which would still be acceptable for many CubeSat missions. In all 100 trials, the 1RW satellite had less than $15^\circ$ of final attitude error.

\section{Analysis}
%\todo{Edit new analysis}
\subsection{Simple 180$^\circ$ Slew}
For the simple $180^\circ$ slews, the median final attitude error for the 1RW satellite in our Monte Carlo simulation is 5.3e-4 degrees, several orders of magnitude lower than the median error of 4.20 degrees for the magnetorquer-only satellite. As we expected, adding a single reaction wheel to a satellite, which allows torque around the magnetic field line in most orientations, can drastically improve pointing performance. 
\\ \\
Although the algorithm found 1RW trajectories with sub-degree pointing precision in only 96\% of all trials, it achieved pointing precision better than 10 degrees in all trials. For most CubeSat missions, this pointing performance is acceptable. For those missions with sub-degree pointing requirements, our extension of ALTRO is fast enough that if a slew is not time-critical, a 1RW satellite could replan and attempt to slew again if the first plan results in more than one degree of final attitude error. Furthermore, all simulated attitude trajectories for the 1RW satellite result in $95\%$ or less of the reaction wheel momentum storage being used -- the reaction wheel never saturates, and in many cases stored angular momentum is driven close to zero by the end of the trajectory.
\\ \\
For the magnetorquer-only satellite, only 73\% of trials in the Monte Carlo simulation resulted in less than 10 degrees of planned final attitude error, and only 11\% of trials resulted in less than one degree of final attitude error. For many CubeSat missions, this higher level of error would necessitate routinely replanning slews in order to guarantee acceptable levels of pointing precision. However, a magnetorquer-only attitude control system could be appropriate for missions with flexible pointing requirements, or with longer time to prepare for slews, execute them, and stabilize after.
\subsection{Reduced Attitude Goal Slews}
Compared to the simple $180^\circ$ slews, the reduced attitude  slews provide a less-restrictive case that is also more realistic for CubeSats -- the case where a satellite camera needs to be aligned with a ground target, but the rotation around the camera axis or vector towards the ground target is unconstrained. This approach has previously been considered for a variety of attitude control questions.\cite{bullo_control_1995,chaturvedi_rigid-body_2011,pong_reduced-attitude_2015,de_angelis_single-axis_2015,McKeen_Method_2022}. By leaving one rotational dimension unconstrained, these slews correspond to a greater number of feasible trajectories and should allow for easier convergence. This should result in lower pointing error than the simple slews, particularly for the magnetorquer-only satellite, where the underactuation about one axis is somewhat compensated for by the goal being reduced by one axis. 
\\ \\
For the reduced attitude slews, the median final attitude error for the 1RW satellite is 2.55e-3 degrees, an order of magnitude lower than the median error of 0.041 degrees for the magnetorquer-only satellite. Furthermore, as in the case of the simple $180^\circ$ slews, the 1RW satellite achieved less than 10 degrees of pointing error in all trials and sub-degree precision in 96\% of trials. The magnetorquer-only satellite achieved less than 10 degrees of pointing error in 75\% of trials and less than one degree of pointing error in 67\% of trials.
\\\\
While the 1RW case is still superior, the marked increase in MTQ-only performance is worth noting. As will be addressed in other forthcoming publications from the authors, this reduced-attitude pointing goal has substantial benefits for magnetorquer-only satellites. One of the most significant issues with MTQ-only control is the underactuation, specifically, the lack of control about a single axis (aligned with the B field). By changing a fully-specified three-dimensional orientation goal to a reduced-attitude single-axis pointing goal, the goal is reduced by a single axis as well, dramatically increasing the set of rotation trajectories that can meet the goal (and making it easier for the solver to find them). This is true even though the non-actuated axis (the B-field) does not usually align with the reduced-attitude axis. 
%\todo{ \\\\keep or delete?\\\\However, while rotation is three-dimensional, it is noncartesian and the missing dimension can be addressed through rotations around the controllable axes. However, the space of working trajectories is hugely increased by redducing the decgree of the pointing goal, allowing for easier finding of these trajectories that compensate for impossible rotations by combining multiple perpendicualr roataions}

\subsection{Multiple Targets}
The multiple target trajectories provide a realistic but difficult case for CubeSats. Although these slews only constrain one axis of the satellite (similar to the case discussed in Reduced Attitude Goals), the multiple targets and shorter initial slew to the first target make these slews more difficult to execute. Further, they capture a case that is currently rarely implemented on CubeSats and other imaging spacecraft--the case where the satellite is switching between targets, trying to meet multiple goals in a single orbit or even a single pass. (In this case, three goals in 500 seconds.) This could allow for imaging a larger area in higher resolution, by taking multiple pictures back-to-back, or meeting multiple goals in a single pass.
\\ \\
For the multiple target trajectories, the median final attitude error for the 1RW satellite is 0.029 degrees and the median error for the magnetorquer-only satellite is 0.22 degrees. The 1RW satellite achieved less than 10 degrees of pointing error in 98\% of cases for the first target and in 99\% of cases for the second and third target. It achieved less than one degree of pointing error in 37\% of cases for the first target, 56\% of cases for the second target, and 91\% of cases for the third target. In comparison, the magnetorquer-only satellite achieved less than 10 degrees of pointing error in 13\% of cases for the first target, 83\% of cases for the second target, and 91\% of cases for the third target, and achieved less than one degree of pointing error in 1\% of cases for the first target, 0\% of cases for the second target, and 77\% of cases for the third target.
\\ \\
Although the magnetorquer-only satellite performs well when pointing at the third target, it performs quite poorly in pointing at the first two imaging targets, demonstrating that a magnetorquer-only attitude control system is likely inappropriate for missions that require multiple target passes. Furthermore, even in its ``good convergence" case, the magnetorquer-only satellite demonstrates attitude drift during the imaging periods, which could potentially smear images. The 1RW satellite, however, generally achieves less than 10 degrees of pointing error, and in over a third of cases achieves sub-degree precision even on the first pointing target, demonstrating that a 1RW system with trajectory planning could execute multiple target trajectories, allowing for multiple images in a single pass or imaging larger areas at a higher resolution.
\subsection{Summary}
Although the 1RW satellite shows clear gains in pointing performance, there are SWaP-C tradeoffs associated with adding a reaction wheel to a satellite. In this work, we simulate a reaction wheel with the properties of a CubeSpace CubeWheel Small: a cost of \$4700, a mass of 60g, dimensions of 23 x 31 x 26 mm, and a peak power consumption of 600 mW \cite{cubewheel}. These properties are typical of a CubeSat-class commercial-off-the-shelf (COTS) reaction wheel. CubeSats with extremely tight budgets for mass, volume, power, or cost may not be able to incorporate such a component. Furthermore, integrating a reaction wheel into a CubeSat bus requires mechanical, electrical, and software engineering work, adding to labor costs and the timeline of a project. As such, the magnetorquer-only attitude control system may be more appropriate for missions with flexible pointing requirements and tight budgets or timelines.
\section{Conclusions \& Future Work}
In this paper, we extend the trajectory planning approach for controlling magnetorquer-only satellites pioneered by Gatherer and Manchester\cite{gatherer2018magnetorquer} to demonstrate a novel technique for controlling small satellites with three magnetorquers and one reaction wheel. We demonstrate clear improvements in pointing performance from adding a single reaction wheel to an otherwise magnetorquer-only system for three types of slews: simple $180^\circ$ slews, reduced attitude goal slews, and multi-target slews. 1RW systems have lower SWaP-C than systems with three reaction wheels and three magnetorquers but have excellent pointing performance, and as such present an attractive low-cost but high-performance attitude control option for CubeSats. Additionally, the use of trajectory planning allows for intelligent application of the reaction wheel and constant management of angular momentum and stiction avoidance, thus reducing the complexity of operations.
\subsection{Future Work}
BeaverCube-2 is a 3U satellite currently being developed by MIT STAR Lab and the Northrop Grumman Corporation which will fly a 1RW system and use the trajectory planning approach outlined in this paper for attitude control. BeaverCube-2's attitude determination and control system (ADCS) also includes an estimated Kalman filter (EKF) with an extra term for tracking external disturbance torques in order to estimate attitude, and an autonomous system onboard that calculates attitude target trajectories to input into ALTRO given a pointing mode (e.g. nadir), respects angular velocity limits, and avoids pointing cameras at the sun.
\\ \\
We have a number of improvements planned for the 1RW system on BeaverCube-2: we plan to use a model-predictive controller (MPC) rather than an LQR controller to track a planned trajectory, and we plan to implement additional cost functions beyond the one outlined in this paper for different pointing modes. We also plan to add further speed optimizations to our extension of ALTRO in order to run it efficiently on BeaverCube-2's flight hardware. These improvements will allow the attitude control approach outlined in this paper to be operationalized on a CubeSat.
\\ \\
The trajectory optimization approach to attitude control for 1RW satellites outlined in this paper also shows promise for desaturating reaction wheels on both 1RW systems and fully actuated control systems with three reaction wheels (3RW) while maintaining attitude. In the future, we plan to demonstrate persistent pointing while driving stored angular momentum to low levels (not zero, due to stiction) for both 1RW and 3RW systems using ALTRO. 
\\ \\
We are currently refactoring, documenting, and testing our code in preparation for public release and plan to open-source it in the near future in concert with a paper or tutorial that thoroughly explains our implementation of ALTRO.

\section{Acknowledgements}
This material is based upon work supported by the
National Science Foundation Graduate Research Fellowship under Grant No. 2141064. This material is also based upon work completed while supported by MIT Portugal:
Portuguese Science and Technology Foundation for work on the
project “AEROS Constellation—Development of a Nanosatellite
Platform as a Precursor of a Future Constellation to Leverage the
Space/Ocean Scientific and Economic Synergies” and by the Northrop Grumman Corporation for work on BeaverCube-2 hardware development.
\bibliographystyle{AAS_publication}   % Number the references.
\bibliography{references}   % Use references.bib to resolve the labels.

@inproceedings{howell2019altro,
  title={Altro: A fast solver for constrained trajectory optimization},
  author={Howell, Taylor A and Jackson, Brian E and Manchester, Zachary},
  booktitle={2019 IEEE/RSJ International Conference on Intelligent Robots and Systems (IROS)},
  pages={7674--7679},
  year={2019},
  organization={IEEE}
 }

@article{babcock2011cubesat,
  title={CubeSat attitude determination via Kalman filtering of magnetometer and solar cell data},
  author={Babcock, Erik},
  year={2011}
}

@article{sola2017quaternion,
  title={Quaternion kinematics for the error-state Kalman filter},
  author={Sola, Joan},
  journal={arXiv preprint arXiv:1711.02508},
  year={2017}
}

@article{gatherer2018magnetorquer,
  title={Magnetorquer-Only Attitude Control of Small Satellites Using Trajectory Optimization},
  author={Gatherer, Andrew and Manchester, Zac},
  journal={Adv. Astronautical Sci. AAS},
  pages={19--927},
  year={2018}
}

@inproceedings{zheng2015improving,
  title={Improving deep neural networks using softplus units},
  author={Zheng, Hao and Yang, Zhanlei and Liu, Wenju and Liang, Jizhong and Li, Yanpeng},
  booktitle={2015 International joint conference on neural networks (IJCNN)},
  pages={1--4},
  year={2015},
  organization={IEEE}
}

@article{chaturvedi_rigid-body_2011,
	title = {Rigid-{Body} {Attitude} {Control}},
	volume = {31},
	issn = {1941-000X},
	doi = {10.1109/MCS.2011.940459},
	number = {3},
	journal = {IEEE Control Systems Magazine},
	author = {Chaturvedi, Nalin A. and Sanyal, Amit K. and McClamroch, N. Harris},
	month = jun,
	year = {2011},
	note = {Conference Name: IEEE Control Systems Magazine},
	pages = {30--51},
}

@article{gatherer_magnetorquer-only_2019,
	title = {{MAGNETORQUER}-{ONLY} {ATTITUDE} {CONTROL} {OF} {SMALL} {SATELLITES} {USING} {TRAJECTORY} {OPTIMIZATION}},
	language = {en},
	journal = {Proceedings of AAS/AIAA Astrodynamics Specialist Conference},
	author = {Gatherer, Andrew and Manchester, Zac},
	month = aug,
	year = {2019},
	pages = {14},
}

@article{jackson_planning_2021,
	title = {Planning {With} {Attitude}},
	issn = {2377-3766, 2377-3774},
	url = {https://ieeexplore.ieee.org/document/9326337/},
	doi = {10.1109/LRA.2021.3052431},
	language = {en},
	urldate = {2021-07-05},
	journal = {IEEE Robotics and Automation Letters},
	author = {Jackson, Brian E. and Tracy, Kevin and Manchester, Zachary},
	year = {2021},
	pages = {1--1},
}

@inproceedings{howell_altro_2019,
	address = {Macau, China},
	title = {{ALTRO}: {A} {Fast} {Solver} for {Constrained} {Trajectory} {Optimization}},
	isbn = {978-1-72814-004-9},
	shorttitle = {{ALTRO}},
	url = {https://ieeexplore.ieee.org/document/8967788/},
	doi = {10.1109/IROS40897.2019.8967788},
	language = {en},
	urldate = {2021-07-05},
	booktitle = {2019 {IEEE}/{RSJ} {International} {Conference} on {Intelligent} {Robots} and {Systems} ({IROS})},
	publisher = {IEEE},
	author = {Howell, Taylor A. and Jackson, Brian E. and Manchester, Zachary},
	month = nov,
	year = {2019},
	pages = {7674--7679},
}

@article{de_angelis_single-axis_2015,
	title = {Single-{Axis} {Pointing} of {Underactuated} {Spacecraft} in the {Presence} of {Path} {Constraints}},
	volume = {38},
	url = {https://doi.org/10.2514/1.G000121},
	doi = {10.2514/1.G000121},
	number = {1},
	urldate = {2021-08-24},
	journal = {Journal of Guidance, Control, and Dynamics},
	author = {de Angelis, Emanuele L. and Giulietti, Fabrizio and Avanzini, Giulio},
	year = {2015},
	note = {Publisher: American Institute of Aeronautics and Astronautics
\_eprint: https://doi.org/10.2514/1.G000121},
	pages = {143--147},
}

@article{frik_attitude_1970,
	title = {Attitude stability of satellites subjected to gravity gradient and aerodynamic torques},
	volume = {8},
	issn = {0001-1452},
	url = {http://arc.aiaa.org/doi/10.2514/3.5990},
	doi = {10.2514/3.5990},
	number = {10},
	urldate = {2021-07-05},
	journal = {AIAA Journal},
	author = {Frik, Martin A.},
	month = oct,
	year = {1970},
	note = {Publisher: American Institute of Aeronautics and Astronautics},
	pages = {1780--1785},
}

@inproceedings{wang_satellite_1998,
	title = {Satellite attitude control using only magnetorquers},
	volume = {1},
	doi = {10.1109/ACC.1998.694663},
	booktitle = {Proceedings of the 1998 {American} {Control} {Conference}. {ACC} ({IEEE} {Cat}. {No}.{98CH36207})},
	author = {Wang, Ping and Shtessel, Y.B.},
	month = jun,
	year = {1998},
	note = {ISSN: 0743-1619},
	pages = {222--226 vol.1},
}

@book{wisniewski_satellite_1997,
	title = {Satellite {Attitude} {Control} {Using} {Only} {Electromagnetic} {Actuation}},
	url = {https://vbn.aau.dk/en/publications/satellite-attitude-control-using-only-electromagnetic-actuation},
	language = {Danish},
	urldate = {2021-03-19},
	publisher = {Aalborg Universitetsforlag},
	author = {Wisniewski, Rafal},
	year = {1997},
}

@book{wisniewski_three-axis_1996,
	title = {Three-{Axis} {Satellite} {Attitude} {Control} {Based} on {Magnetic} {Torquing} - {Linear} {Optimal} {Approach}},
	author = {Wisniewski, Rafal},
	year = {1996},
}

@misc{m_campbell_beavercube_2019,
	title = {{BeaverCube} {CDR} {CAD}},
	urldate = {2021-11-24},
	author = {{M. Campbell} and {T. Tran}},
	month = dec,
	year = {2019},
}

@misc{jackson_al_ilqr_tutorialpdf_2019,
	title = {{AL}\_iLQR\_Tutorial.pdf},
	url = {https://bjack205.github.io/papers/AL_iLQR_Tutorial.pdf},
	urldate = {2021-11-04},
	journal = {AL-iLQR Tutorial},
	author = {Jackson, Brian},
	month = jun,
	year = {2019},
}

@article{bullo_control_1995,
	series = {3rd {IFAC} {Symposium} on {Nonlinear} {Control} {Systems} {Design} 1995, {Tahoe} {City}, {CA}, {USA}, 25-28 {June} 1995},
	title = {Control on the {Sphere} and {Reduced} {Attitude} {Stabilization}},
	volume = {28},
	issn = {1474-6670},
	url = {https://www.sciencedirect.com/science/article/pii/S1474667017468789},
	doi = {10.1016/S1474-6670(17)46878-9},
	language = {en},
	number = {14},
	urldate = {2021-11-04},
	journal = {IFAC Proceedings Volumes},
	author = {Bullo, F. and Murray, R. M. and Sarti, A.},
	month = jun,
	year = {1995},
	pages = {495--501},
}

@article{wisniewski_linear_2000,
	title = {Linear {Time}-{Varying} {Approach} to {Satellite} {Attitude} {Control} {Using} {Only} {Electromagnetic} {Actuation}},
	volume = {23},
	url = {http://arc.aiaa.org/doi/10.2514/2.4609},
	doi = {10.2514/2.4609},
	number = {4},
	urldate = {2021-11-03},
	journal = {Journal of Guidance, Control, and Dynamics},
	author = {Wisniewski, Rafal},
	month = jul,
	year = {2000},
	note = {Publisher: American Institute of Aeronautics and Astronautics},
	pages = {640--647},
}

@article{psiaki_magnetic_2001,
	title = {Magnetic {Torquer} {Attitude} {Control} via {Asymptotic} {Periodic} {Linear} {Quadratic} {Regulation}},
	volume = {24},
	url = {http://arc.aiaa.org/doi/10.2514/2.4723},
	doi = {10.2514/2.4723},
	number = {2},
	urldate = {2021-11-03},
	journal = {Journal of Guidance, Control, and Dynamics},
	author = {Psiaki, Mark L.},
	month = mar,
	year = {2001},
	note = {Publisher: American Institute of Aeronautics and Astronautics},
	pages = {386--394},
}

@article{wisniewski_sliding_1998,
	series = {14th {IFAC} {Symposium} on {Automatic} {Control} in {Aerospace} 1998, {Seoul}, {Korea}, 24-28 {August} 1998},
	title = {Sliding {Mode} {Attitude} {Control} for {Magnetic} {Actuated} {Satellite}},
	volume = {31},
	issn = {1474-6670},
	url = {https://www.sciencedirect.com/science/article/pii/S1474667017410767},
	doi = {10.1016/S1474-6670(17)41076-7},
	language = {fr},
	number = {21},
	urldate = {2021-11-03},
	journal = {IFAC Proceedings Volumes},
	author = {Wiśniewski, Rafał},
	month = aug,
	year = {1998},
	pages = {179--184},
}

@article{lovera_global_2005,
	title = {Global {Magnetic} {Attitude} {Control} of {Inertially} {Pointing} {Spacecraft}},
	volume = {28},
	url = {http://arc.aiaa.org/doi/10.2514/1.11844},
	doi = {10.2514/1.11844},
	number = {5},
	urldate = {2021-11-03},
	journal = {Journal of Guidance, Control, and Dynamics},
	author = {Lovera, Marco and Astolfi, Alessandro},
	month = sep,
	year = {2005},
	note = {Publisher: American Institute of Aeronautics and Astronautics},
	pages = {1065--1072},
}

@article{belokonov_analysis_2018,
	title = {Analysis and {Synthesis} of {Motion} of {Aerodynamically} {Stabilized} {Nanosatellites} of the {CubeSat} {Design}},
	volume = {9},
	issn = {2075-1109},
	url = {https://doi.org/10.1134/S2075108718040028},
	doi = {10.1134/S2075108718040028},
	language = {en},
	number = {4},
	urldate = {2021-11-03},
	journal = {Gyroscopy and Navigation},
	author = {Belokonov, I. V. and Timbai, I. A. and Nikolaev, P. N.},
	month = oct,
	year = {2018},
	pages = {287--300},
}

@inproceedings{stuurman_ryefemsat_2010,
	series = {{SpaceOps} {Conferences}},
	title = {{RyeFemSat}: {Ryerson} {University} {Femtosatellite} {Design}  and {Testing}},
	shorttitle = {{RyeFemSat}},
	url = {http://arc.aiaa.org/doi/10.2514/6.2010-2157},
	doi = {10.2514/6.2010-2157},
	urldate = {2021-11-03},
	booktitle = {{SpaceOps} 2010 {Conference}},
	publisher = {American Institute of Aeronautics and Astronautics},
	author = {Stuurman, Bryan and Kumar, Krishna},
	month = apr,
	year = {2010},
}

@article{lovera_global_2006,
	title = {Global magnetic attitude control of spacecraft in the presence of gravity gradient},
	volume = {42},
	issn = {1557-9603},
	doi = {10.1109/TAES.2006.248214},
	number = {3},
	journal = {IEEE Transactions on Aerospace and Electronic Systems},
	author = {Lovera, Marco and Astolfi, Alessandro},
	month = jul,
	year = {2006},
	note = {Conference Name: IEEE Transactions on Aerospace and Electronic Systems},
	pages = {796--805},
}

@inproceedings{chasset_3-axis_2006,
	series = {International {Astronautical} {Congress} ({IAF})},
	title = {3-axis magnetic control with multiple attitude profile capabilities in the {PRISMA} mission},
	url = {http://arc.aiaa.org/doi/10.2514/6.IAC-06-C1.2.03},
	doi = {10.2514/6.IAC-06-C1.2.03},
	urldate = {2021-11-03},
	booktitle = {57th {International} {Astronautical} {Congress}},
	publisher = {American Institute of Aeronautics and Astronautics},
	author = {Chasset, Camille and Berge, Sten and Bodin, Per and Jakobsson, Björn},
	month = oct,
	year = {2006},
}

@book{cappelletti_cubesat_2020,
	title = {{CubeSat} {Handbook}: {From} {Mission} {Design} to {Operations}},
	isbn = {978-0-12-817884-3},
	shorttitle = {{CubeSat} {Handbook}},
	author = {Cappelletti, Chantal and Battistini, Simone and Malphrus, Benjamin},
	month = sep,
	year = {2020},
}

@book{wertz_space_2011,
	address = {Hawthorne, CA},
	title = {Space mission engineering : the new {SMAD}},
	isbn = {1-881883-16-7 978-1-881883-16-6 1-881883-15-9 978-1-881883-15-9},
	language = {English},
	publisher = {Microcosm Press : Sold and distributed worldwide by Microcosm Astronautics Books},
	author = {Wertz, Jame Richard and Everett, David F. and Puschell, Jeffery John},
	year = {2011},
}

@article{ovchinnikov_survey_2019,
	title = {A survey on active magnetic attitude control algorithms for small satellites},
	volume = {109},
	doi = {10.1016/j.paerosci.2019.05.006},
	journal = {Progress in Aerospace Sciences},
	author = {Ovchinnikov, Michael and Roldugin, Dmitry},
	month = jun,
	year = {2019},
}

@article{wisniewski_optimal_1999,
	series = {14th {IFAC} {World} {Congress} 1999, {Beijing}, {Chia}, 5-9 {July}},
	title = {Optimal {Magnetic} {Attitude} {Control}},
	volume = {32},
	issn = {1474-6670},
	url = {https://www.sciencedirect.com/science/article/pii/S1474667017573632},
	doi = {10.1016/S1474-6670(17)57363-2},
	language = {en},
	number = {2},
	urldate = {2021-11-03},
	journal = {IFAC Proceedings Volumes},
	author = {Wisniewski, Rafal and Landis Markley, F.},
	month = jul,
	year = {1999},
	pages = {7991--7996},
}

@article{chasset_3-axis_2013,
	title = {3-{Axis} magnetic control: flight results of the {TANGO} satellite in the {PRISMA} mission},
	volume = {5},
	issn = {1868-2510},
	shorttitle = {3-{Axis} magnetic control},
	url = {https://doi.org/10.1007/s12567-013-0034-9},
	doi = {10.1007/s12567-013-0034-9},
	language = {en},
	number = {1},
	urldate = {2021-11-02},
	journal = {CEAS Space Journal},
	author = {Chasset, C. and Noteborn, R. and Bodin, P. and Larsson, R. and Jakobsson, B.},
	month = sep,
	year = {2013},
	pages = {1--17},
}

@article{ivanov_advanced_2017,
	title = {Advanced numerical study of the three-axis magnetic attitude control and determination with uncertainties},
	volume = {132},
	issn = {0094-5765},
	url = {https://www.sciencedirect.com/science/article/pii/S009457651630371X},
	doi = {10.1016/j.actaastro.2016.11.045},
	language = {en},
	urldate = {2021-11-01},
	journal = {Acta Astronautica},
	author = {Ivanov, D. S. and Ovchinnikov, M. Yu. and Penkov, V. I. and Roldugin, D. S. and Doronin, D. M. and Ovchinnikov, A. V.},
	month = mar,
	year = {2017},
	pages = {103--110},
}

@article{sechi_magnetic_2005,
	title = {Magnetic {Attitude} {Control} of the {Goce} {Satellite}},
	volume = {606},
	author = {Sechi, G. and André, G. and Andreis, D. and Saponara, M.},
	month = dec,
	year = {2005},
	pages = {37},
}

@inproceedings{krogstad_explicit_2005,
	title = {Explicit {Model} {Predictive} {Control} of a {Satellite} with {Magnetic} {Torquers}},
	doi = {10.1109/.2005.1467064},
	booktitle = {Proceedings of the 2005 {IEEE} {International} {Symposium} on, {Mediterrean} {Conference} on {Control} and {Automation} {Intelligent} {Control}, 2005.},
	author = {Krogstad, T. and Gravdahl, J.T. and Tondel, P.},
	month = jun,
	year = {2005},
	note = {ISSN: 2158-9879},
	pages = {491--496},
}

@inproceedings{liang_time-optimal_2004,
	title = {Time-optimal magnetic attitude control for small spacecraft},
	volume = {1},
	doi = {10.1109/CDC.2004.1428639},
	booktitle = {2004 43rd {IEEE} {Conference} on {Decision} and {Control} ({CDC}) ({IEEE} {Cat}. {No}.{04CH37601})},
	author = {Liang, Jinsong and Fullmer, R. and Chen, Yang Quan},
	month = dec,
	year = {2004},
	note = {ISSN: 0191-2216},
	pages = {255--260 Vol.1},
}

@article{kumar_h_2016,
	series = {4th {IFAC} {Conference} on {Advances} in {Control} and {Optimization} of {Dynamical} {Systems} {ACODS} 2016},
	title = {H$\inf$ {Tracking} {Control} for {Magnetically} {Controlled} {Nano}-satellite},
	volume = {49},
	issn = {2405-8963},
	url = {https://www.sciencedirect.com/science/article/pii/S2405896316300477},
	doi = {10.1016/j.ifacol.2016.03.047},
	language = {en},
	number = {1},
	urldate = {2021-10-20},
	journal = {IFAC-PapersOnLine},
	author = {Kumar, Shilpee and Varma, Mithun S. and Rao, A. Divya and Agrawal, V. K.},
	month = jan,
	year = {2016},
	pages = {166--172},
}

@inproceedings{janardhanan_attitude_2012,
	title = {Attitude control of magnetic actuated spacecraft using super-twisting algorithm with nonlinear sliding surface},
	doi = {10.1109/VSS.2012.6163476},
	booktitle = {2012 12th {International} {Workshop} on {Variable} {Structure} {Systems}},
	author = {Janardhanan, S and un Nabi, Mashuq and Tiwari, Pyare Mohan},
	month = jan,
	year = {2012},
	note = {ISSN: 2158-3986},
	pages = {46--51},
}

@inproceedings{sofyali_purely_2012,
	title = {Purely magnetic spacecraft attitude control by using classical and modified sliding mode algorithms},
	doi = {10.1109/VSS.2012.6163488},
	booktitle = {2012 12th {International} {Workshop} on {Variable} {Structure} {Systems}},
	author = {Sofyali, Ahmet and Jafarov, Elbrous M.},
	month = jan,
	year = {2012},
	note = {ISSN: 2158-3986},
	pages = {117--123},
}

@article{pong_reduced-attitude_2015,
	title = {Reduced-{Attitude} {Boresight} {Guidance} and {Control} on {Spacecraft} for {Pointing}, {Tracking}, and {Searching}},
	volume = {38},
	issn = {0731-5090, 1533-3884},
	url = {https://arc.aiaa.org/doi/10.2514/1.G000264},
	doi = {10.2514/1.G000264},
	language = {en},
	number = {6},
	urldate = {2021-10-13},
	journal = {Journal of Guidance, Control, and Dynamics},
	author = {Pong, Christopher M. and Miller, David W.},
	month = jun,
	year = {2015},
	pages = {1027--1035},
}

@article{McKeen_Method_2022,
author = {McKeen, Patrick and Cahoy, Kerri},
title = {Method for Pointing Spacecraft While Minimizing Change in Orientation},
journal = {Journal of Guidance, Control, and Dynamics},
volume = {45},
number = {8},
pages = {1554-1557},
year = {2022},
doi = {10.2514/1.G006429},

URL = { 
        https://doi.org/10.2514/1.G006429
    
},
eprint = { 
        https://doi.org/10.2514/1.G006429
    
}

}

@manual{cubewheel,
    organization  = "CubeSpace",
    title         = "CubeWheel specification",
    url = {https://www.cubespace.co.za/products/adcs-components/cubewheel/}
}

\end{document}